\documentclass[review]{elsarticle}

\usepackage{lineno,hyperref}

\usepackage{pstricks}
\usepackage{pstricks-add}
\usepackage{pst-plot}
\usepackage{amssymb}
\usepackage{graphicx}
\usepackage{caption}
\usepackage{subcaption}
\usepackage{verbatim}
\usepackage[english]{biblio}
\usepackage{enumerate}
\usepackage{amsmath}
\usepackage{psfrag}
\usepackage{tikz}
\usepackage{verbatim}
\newtheorem{thm}{Theorem}[section]

\newtheorem{rem}[thm]{Remark}

\newcommand{\Rr}{\mathbb R}

\newcommand{\bu}{{\bf u}}
\newcommand{\bv}{{\bf v}}

\newcommand{\be}{\begin{equation}}
\newcommand{\ee}{\end{equation}}

\def\Xint#1{\mathchoice
{\XXint\displaystyle\textstyle{#1}}%
{\XXint\textstyle\scriptstyle{#1}}%
{\XXint\scriptstyle\scriptscriptstyle{#1}}%
{\XXint\scriptscriptstyle\scriptscriptstyle{#1}}%
\!\int}
\def\XXint#1#2#3{{\setbox0=\hbox{$#1{#2#3}{\int}$}
\vcenter{\hbox{$#2#3$}}\kern-.5\wd0}}

\def\dashint{\Xint-}

\usepackage{hyperref}
\hypersetup{
  colorlinks   = true, 
  urlcolor     = blue, 
  linkcolor    = blue, 
  citecolor   = purple 
}

\numberwithin{equation}{section}

\begin{document}

\begin{frontmatter}

\title{High-order In-cell Discontinuous Reconstruction path-conservative methods for nonconservative hyperbolic systems - DR.MOOD generalization}


\author[uma]{Ernesto Pimentel-Garc\'ia\corref{mycorrespondingauthor}}
\cortext[mycorrespondingauthor]{Corresponding author}
\ead{erpigar@uma.es}

\author[uma]{Manuel J. Castro}
\ead{mjcastro@uma.es}

\author[versailles]{Christophe Chalons}
\ead{christophe.chalons@uvsq.fr}


\author[uma]{Carlos Par\'es}
\ead{pares@uma.es}

\address[uma]{Departamento  de  An\'alisis  Matem\'atico,  Estad\'istica  e  Investigaci\'on  Operativa,  y  Matem\'atica  aplicada,  Universidad  de  M\'alaga,Bulevar Louis Pasteur, 31, 29010, M\'alaga, Spain. }
\address[versailles]{Laboratoire de Math\'ematiques de Versailles, UVSQ, CNRS, Universit\'e Paris-Saclay, 78035 Versailles, France.}

\begin{abstract}
In this work we develop a new framework to deal numerically with discontinuous solutions in nonconservative hyperbolic systems. First an extension of the MOOD methodology to nonconservative systems based on Taylor expansions is presented. This extension combined with an in-cell discontinuous reconstruction operator are the key points to develop a new family of high-order methods that are able to capture exactly isolated shocks. Several test cases are proposed to validate these methods for the Modified Shallow Water equations and the Two-Layer Shallow Water system.
\end{abstract}


\end{frontmatter}


\tableofcontents

\section{Introduction}\label{S:intro}

  We consider  first
order quasi-linear PDE systems
\begin{equation}\label{sys:nonconservative}
\partial_t {\bf u} + {\mathcal{A}}({\bf u}) \partial_x {\bf u} = 0, \quad x \in \Rr, \quad t \in \Rr^{+},
\end{equation}
in which the unknown $\bu(x,t)$ takes values in an open convex set
$\Omega$ of $\mathbb{R}^N$, and
$\mathcal{A}(\bf u)$ is a smooth locally bounded map from $\Omega$ to $\mathcal{M}_{N\times
N}(\mathbb{R})$. We focus on strictly hyperbolic systems and the characteristic fields
$R_i(\bu)$, $i=1,\ldots ,N$, are supposed to be either genuinely
nonlinear:
$$
\nabla \lambda _i(\bu) \cdot R_i(\bu) \neq 0, \quad \forall \,
\bu\in\Omega,
$$
or linearly degenerate:
$$
\nabla \lambda _i(\bu) \cdot R_i(\bu) = 0, \quad \forall \, \bu \in \Omega.
$$
Here, $\lambda_1(\bu)\leq\ldots\leq\lambda_N(\bu)$ represent the eigenvalues
of $\mathcal{A}(\bu)$ and $R_1(\bu), \ldots,
R_N(\bu)$ a set of associated eigenvectors.

Different mathematical theories allow one to define the nonconservative product $\mathcal{A}(\bu)\,\partial_x \bu$. Here the theory developed by Dal Maso, LeFloch, and Murat
\cite{Maso1995} that allows a notion of weak solution which satisfies \eqref{sys:nonconservative} in the sense of Borel measures is followed. This definition is based on a family of Lipschitz continuous paths $\Phi: [0,1]
\times \Omega \times \Omega \to \Omega$, which must
satisfy certain regularity and  compatibility conditions, in particular
\begin{equation} \label{cond1}
\Phi(0;\bu_l, \bu_r) = \bu_l, \qquad  \Phi(1;\bu_l, \bu_r) = \bu_r,\quad \forall \bu_l, \bu_r \in \Omega,
\end{equation}
and
\begin{equation} \label{cond2}
\Phi(s; \bu, \bu) = \bu, \quad \forall \bu \in \Omega.
\end{equation}
The interested reader
is addressed to \cite{Maso1995}  for a rigorous and complete
presentation of this theory. 
 The family of paths can be understood as a tool to
give a sense to integrals of the form
$$
\int_a^b \mathcal{A}(\bu(x))\,\partial_x \bu(x) \, dx,
$$
for functions $\bu$ with jump discontinuities. More precisely,
given a bounded variation function $\bu:[a,b] \mapsto \Omega$, we define:
\begin{equation}\label{wint}
\dashint_a^b \mathcal{A}(\bu(x))\,\partial_x \bu(x) \, dx =   \int_a^b \mathcal{A}(\bu(x))\,\partial_x \bu(x) \,dx  
+ \sum_m \int_0^1\mathcal{A}(\Phi(s;\bu_m^-,\bu_m^+))\frac{\partial\Phi}{\partial
s}(s;\bu_m^-,\bu_m^+)\,ds,
\end{equation}
where $\bu_m^-$ and $\bu_m^+$ represent, respectively,  the
limits of $\bu$ to the left and right of its $m$th discontinuity. Observe that, in \eqref{wint}, the family
of paths has been used to determine the Dirac measures placed at
the discontinuities of $\bu$.

If such a mathematical definition of the nonconservative products is assumed to define the concept of weak solution,  the generalized Rankine-Hugoniot condition:
\begin{equation} \label{gen_R-H}
 \int_0^1 \mathcal{A}(\Phi(s;\bu^-\bu^+))
\frac{\partial\Phi}{\partial s}(s;\bu^-,\bu^+)\, ds  = \sigma (\bu^+ - \bu^-)
\end{equation}
has to be satisfied across an admissible discontinuity. Here,
$\sigma $ is the speed of propagation of the discontinuity, and
$\bu^-$ and $\bu^+$ are the left and right limits of the solution at
the discontinuity. 

Once the family of paths has been prescribed, a concept of entropy
is required, as it happens for systems of conservation laws, that may be given by an entropy pair  or by Lax's entropy criterion. 

As we observe in \eqref{gen_R-H} the concept of weak solution depends clearly on the choice of the family
of paths, which is a priori arbitrary, so that the crucial question is
how to choose the 'good' family of paths. In fact, when the
hyperbolic system is the vanishing-viscosity limit of the
parabolic problems
\begin{equation}
\partial_t \bu^\eps + \mathcal{A}(\bu^\eps)\,\partial_x \bu^\eps = \eps
(\mathcal{R}(\bu^\eps)\partial_x \bu^\eps)_x, \label{NC-p}
\end{equation}
where $\mathcal{R}(\bu)$ is any positive-definite matrix, 
the adequate family of paths should be related to the \textit{viscous
profiles}: a function $\bv$ is said to be
a viscous profile for \eqref{NC-p} linking the states $\bu^-$ and  $\bu^+$
if it satisfies
\begin{equation} \label{vpbc}
\lim_{\chi \to -\infty} \bv(\chi ) = \bu^-, \quad  \lim_{\chi \to
+\infty} \bv(\chi ) = \bu^+,  \quad \lim_{\chi \to \pm \infty} \bv'(\chi) = 0
\end{equation}
and there exists $\sigma \in \mathbb{R}$ such that the travelling wave
\begin{equation}
\label{vprofile} \bu^\eps(x,t) = \bv\left( \frac{x - \sigma
t}{\eps}\right),
\end{equation}
is a solution of \eqref{NC-p} for every $\epsilon$.
It can be easily
verified that, in order to be a viscous profile, $\bv$ has to solve the equation
\begin{equation} \label{NC-vpec}
-\sigma \bv' + \mathcal{A}(\bv)\,\bv' = (\mathcal{R}(\bv)\,\bv')',
\end{equation}
with boundary conditions \eqref{vpbc}.
If there exists a viscous profile linking the states $\bu^-$ and $\bu^+$, the good choice for the path connecting the
states would be, after a reparameterization, the
viscous profile $\bv$. 

Unlike standard conservative systems, where the usual Rankine-Hugoniot conditions are always obtained regardless of the viscous term chosen, nonconservative systems exhibit varying jump conditions depending on the selected viscous term $\mathcal{R}$.

Defining numerical methods that accurately converge to the correct weak solutions for system \eqref{sys:nonconservative} is not a simple task. While Lax's equivalence theorem guarantees the right convergence of numerical methods when consistency and stability are fulfilled for linear systems, this does not hold in general for nonlinear problems. For instance, in the case of systems of conservation laws,  even stable conservative methods can converge to solutions that are not admissible weak solution: this is the case for Roe method that may converge to weak solutions that are not entropy solutions. To ensure convergence to the right weak solutions, beyond consistency and stability, entropy has to be well controlled: for instance, entropy-fix techniques have to be added to Roe method (see for example \cite{Harte83}). However, in the case of nonconservative systems, ensuring convergence to the correct weak solutions requires more than consistency, stability, and entropy control. It is also crucial to properly control numerical viscosity and, more generally, the effects of numerical dissipation. Managing these aspects is essential to achieve the correct convergence (see \cite{leflochmishra} for a review on this topic).

Developing finite-difference or finite-volume methods that satisfy all four properties mentioned earlier is challenging in general. The framework of path-conservative methods introduced in \cite{pares2006numerical} makes easy the design of numerical methods that are consistent with the definition of weak-solutions corresponding to a selected family of paths. Moreover, it allows to extend to nonconservative systems well-known families of conservative methods preserving their stability and entropy properties. The numerical solutions provided by such methods converge to limits that are entropy weak solutions but,  if the numerical dissipation is not controlled, they can be weak solutions according to a family of paths different from the selected one. This means that the limits of the numerical solutions are classical solutions in smooth regions but exhibit discontinuities that satisfy a different jump condition \eqref{gen_R-H} from the expected one: see \cite{Castro2008}, \cite{Abgrall2010}. In fact, the family of paths that governs the jump conditions satisfied by the limits of the numerical solutions is related to the viscous profiles of the equivalent equation of the method, as discussed in \cite{Castro2008}. If, for instance, the family of paths is based on the viscous profiles related to a regularization \eqref{NC-p}, the leading terms in the equivalent equation that represent the numerical viscosity of the scheme may not match the viscous term in \eqref{NC-p}.

Different techniques have been introduced to address this convergence issue, at least partially: \cite{Berthon},
\cite{BLFMP}, \cite{audebert2006hybrid}, \cite{berthon2002nonlinear}, \cite{castro2013entropy}, \cite{chalons2005navier}, \cite{chalons2017new}, \cite{fjordholm2012accurate},  \cite{chalons2019path}.  Among these techniques, the path-conservative entropy stable methods presented in \cite{castro2013entropy} and their extension to high-order Discontinuous Galerkin (DG) methods in \cite{hiltebrand} have significantly reduced convergence errors: to do this, entropy-conservative numerical methods are first introduced so the stabilization is achieved through the discretization of the viscous term in the regularized equation \eqref{NC-p}.  


Recently, in \cite{chalons2019path}, an in-cell discontinuous reconstruction technique has been added to first-order path-conservative methods for nonconservative systems. This technique allows one to  capture correctly weak solutions with isolated shock waves. It is based on the discontinuous reconstruction operators introduced for conservative systems in \cite{despres2001contact}, \cite{lagoutiere2004non},\cite{lagoutiere2008stability}. In \cite{pimentel2021cell} this technique was combined with the MUSCL-Hancock strategy to obtain second-oder methods that correctly capture isolated shock waves. The goal of this article is to extend these  methods to more complex systems and to higher-order accuracy. To do this, we combine unlimited  high-order explicit path-conservative numerical methods with the Multi-dimensional Optimal Order Detection (MOOD) strategy (see \cite{diot2012improved, clain2011high}) using an in-cell discontinuous reconstruction method as the parachute robust first-order scheme. In the MOOD strategy, the unlimited high-order explicit scheme is used to produce a candidate numerical solution and the validity of the numerical approximation at every cell is checked using some detection criteria. While the valid ones are kept, the invalid ones are marked. The in-cell discontinuous reconstruction first-order method is then used to recompute the numerical solution in the marked cells. The detection criteria are chosen so that all the cells in which the solution is discontinuous are marked: we obtain therefore high-order accuracy in smooth regions while capturing correctly the discontinuities.  Up to our knowledge, this is the first time that MOOD is used in this framework.

The paper is organized as follows: in Section \ref{sec:O1_pc} a brief introduction to path-conservative methods is given, then in Section \ref{h_o_prediction_step} we introduce high-order path-conservative methods based on Taylor expansions. In Section \ref{sec:mood} the MOOD methodology well adapted to ensure the convergence towards the correct weak solution is described, , giving a detailed description of the correction terms added to be path-conservative. The first-order in-cell discontinuous reconstruction method introduced in \cite{chalons2019path} is recalled in Section \ref{sec:O1_DR} where some modifications are also introduced. Methods that combine the three ingredients (unlimited high-order-path conservative methods, first-order in-cell discontinuous reconstruction schemes, and MOOD limiting), that are called DR.MOOD methods, are introduced in Section \ref{sec:h-o_DRMOOD}. Finally, in Section \ref{sec:numerical_tests} some numerical tests are shown to check the convergence of the methods to the right solutions. In order to do this we consider a nonconservative modified shallow water system and the two-layer shallow water system.



\section{MOOD methods for nonconservative systems.}\label{sec:mood}

The MOOD strategy is based on a hierarchy of methods going from a high-order method with maximal order of accuracy to a robust first-order scheme passing through several methods of decreasing accuracy. In this section, the application of this strategy to nonconservative systems is described to achieve the convergence towards the correct weak solutions. One of the main difficulty when we consider MOOD strategy in nonconservative hyperbolic systems  comes from the numerical treatment of the cells placed at the boundaries of two regions in which numerical methods of different accuracy are used to update the numerical solution. Here we only consider for simplicity, the MOOD strategy based on only two methods: an unlimited high-order method used as a predictor and a first-order scheme used as the corrector in marked cells. Besides these methods, a criterion to mark wrong cells after the prediction step is required as well as the already mentioned correction at the boundary cells.

\subsection{First-order path-conservative methods}\label{sec:O1_pc}
Following \cite{pares2006numerical}, a numerical method for solving \eqref{sys:nonconservative} is said to be path-conservative if  it can be written under the form
\begin{equation}\label{scheme}
\bu_j^{n+1}=\bu_j^n-\frac{\Delta t}{\Delta x}\big(\mathcal{D}_{j-1/2}^{+} +
\mathcal{D}_{j+1/2}^{-}\big),
\end{equation}
where the following notation is used: 
\begin{itemize}

\item $\Delta x$ and $\Delta t$ are the space and time steps respectively. They are supposed to be constant for simplicity. 

\item  $I_j=\left[x_{j-\frac{1}{2}},x_{j+\frac{1}{2}}\right]$ are the computational cells, whose length  is $\Delta x$.

\item $t^n = n \Delta t $, $n = 0, 1 \dots$.

\item $\bu_j^n$ is the approximation of the average of the exact solution at the $j$th cell at time $t^n$, that is, 
\be
\bu_j^n \approx \frac{1}{\Delta x} \int_{x_{j-\frac{1}{2}}}^{x_{j+\frac{1}{2}}}  \bu(x,t^n) \, dr. 
\ee
\item Finally,
$$\mathcal{D}_{j+1/2}^{\pm}=\mathcal{D}^{\pm}\big(\bu_{j}^n,\bu_{j+1}^n \big), $$
where
$\mathcal{D}^-$ and $\mathcal{D}^+$ two Lipschitz continuous functions from $\Omega \times \Omega$ to $\Omega$ that satisfy
\begin{equation}\label{pc1}
\mathcal{D}^{\pm}(\bu,\bu)=0, \quad \forall \bu\in \Omega,
\end{equation}
and
\begin{equation}\label{pc2}
\mathcal{D}^-(\bu_l,\bu_r) + \mathcal{D}^+(\bu_l, \bu_r) = \int _0^1 \mathcal{A}\big(\Phi (s;\bu_l,\bu_r)\big)\frac{\partial
\Phi }{\partial s} (s;\bu_l,\bu_r)\,ds,
\end{equation}
for every  set $\bu_{l}, \bu_r \in \Omega$. 

\end{itemize}

The definition of path-conservative methods is a \textit{formal concept of consistency} for weak solutions defined on the basis of the family of paths $\Phi$: see \cite{castro2017well} for a recent review. In fact, this is a natural extension of the definition of conservative methods for systems of conservation laws: if  $\mathcal{A}(\bu)$ is the Jacobian of a flux function ${F}(\bu)$ it can be easily checked that \eqref{scheme} is equivalent to the conservative method corresponding to the consistent numerical flux given by
\begin{equation}\label{numflux1}
\mathcal{F}(\bu_l, \bu_r) = {F}(\bu_r) - \mathcal{D}^+(\bu_l, \bu_r),
\end{equation}
or, equivalently
\begin{equation}\label{numflux2}
\mathcal{F}(\bu_l, \bu_r) = {F}(\bu_l) + \mathcal{D}^-(\bu_l, \bu_r).
\end{equation}

This framework makes it easy to extend many well-known conservative schemes to nonconservative systems . Let us show three examples:

\begin{itemize}


\item  Godunov method: 
\begin{eqnarray} \label{NCFG1}
\mathcal{D}^-_{G} (\bu_l, \bu_r) =  \int_0^1\mathcal{A}(\Phi(s;\bu_l,\bu_0))\frac{\partial\Phi}{\partial
s}(s;\bu_l,\bu_0)\,ds,\\
\label{NCFG2}
\mathcal{D}^+_{G} (\bu_l, \bu_r) =  \int_0^1\mathcal{A}(\Phi(s;\bu_0,\bu_r))\frac{\partial\Phi}{\partial
s}(s;\bu_0,\bu_r)\,ds,
\end{eqnarray}
where $\bu_0$ is the value at $x=0$ of the
self-similar solution of the Riemann problem 
\begin{equation}\label{RiemannPb}
\left\{
\begin{array}{l}
\displaystyle \partial_t {\bf u} + {\mathcal{A}}({\bf u}) \partial_x {\bf u} = 0, \\
\bu(x, 0) = \begin{cases}
\bu_l &\text{if $x < 0$,}\\
\bu_r & \text{otherwise.}
\end{cases}
\end{array}
\right.
\end{equation}
If the family of paths satisfies some conditions of
compatibility with the solutions of the Riemann problems, the method can be interpreted in terms of the averages of the exact solutions of local Riemann problems in the cells,
as it happens for system of conservation laws:  see
\cite{munoz2007godunov}. 

\item Roe methods:
\begin{equation}\label{NCFR} \mathcal{D}_R^\pm (\bu_l, \bu_r) =
\mathcal{A}^\pm_\Phi(\bu_l, \bu_r) \, (\bu_r - \bu_l),
\end{equation}
where  $\mathcal{A}_\Phi(\bu_l, \bu_r)$ is a  Roe linearization of $\mathcal{A}(\bu)$  in
the sense defined by Toumi in \cite{toumi1992weak}, i.e. a function
$\mathcal{A}_\Phi\colon\Omega \times\Omega \mapsto\mathcal{M}_{N\times
N}(\mathbb{R})$ satisfying the following properties:
\begin{itemize}
\item  for each $\bu_l,\bu_r\in\Omega$,
$\mathcal{A}_\Phi(\bu_l,\bu_r)$ has $N$ distinct real eigenvalues
$\lambda_1(\bu_l, \bu_r)$, \dots, $\lambda_N(\bu_l, \bu_r)$;
\item $\mathcal{A}_\Phi(\bu,\bu)=\mathcal{A}(\bu)$, for every
$\bu\in\Omega$;
\item  for any $\bu_l,\bu_r\in\Omega$,
\begin{equation}\label{Roe-gen}
\mathcal{A}_\Phi(\bu_l,\bu_r)\,
(\bu_r-\bu_l)=\int_0^1\mathcal{A}(\Phi(s;\bu_l,\bu_r))\frac{\partial\Phi}{\partial
s}(s;\bu_l,\bu_r)\,ds.
\end{equation}
\end{itemize}
As usual $\mathcal{A}^\pm_\Phi(\bu_l, \bu_r)$ represent  the matrices whose
eigenvalues are the positive/negative parts of $\lambda_1(\bu_l, \bu_r)$, \dots, $\lambda_N(\bu_l, \bu_r)$
with same eigenvectors. 

\item Rusanov methods:
\begin{equation}\label{NCFRus} \mathcal{D}_{Rus}^\pm (\bu_l, \bu_r) =\left( \mathcal{A}_\Phi(\bu_l,\bu_r)
\pm \frac{1}{2} \max_{k}|\lambda_{k}(\bu_{l}, \bu_{r})| \mathcal{I}\right)(\bu_{r}-\bu_{l}),
\end{equation}
where $\lambda_{k}(\bu_{l}, \bu_{r})$, $k=1,...,N$, are the eigenvalues of the roe matrix $\mathcal{A}_{\Phi}(\bu_{l}, \bu_{r})$ and $\mathcal{I}$ is the identity matrix.
\end{itemize}

\subsection{Unlimited high-order path-conservative methods}\label{h_o_prediction_step}

First-order path-conservative numerical schemes can be extended to high-order by using reconstruction operators:
\begin{equation}\label{eq:hopc}
    {\bf u}_{j}'(t) = -\frac{1}{\Delta x}\left(\mathcal{D}_{j+\frac{1}{2}}^{-}(t) + \mathcal{D}_{j-\frac{1}{2}}^{+}(t) + \int_{x_{j-\frac{1}{2}}}^{x_{j+\frac{1}{2}}}\mathcal{A}(P^t_{j}(x))\frac{\partial}{\partial x}P^t_{j}(x)\, dx\right),
\end{equation}
where:
\begin{itemize}
	\item $\bu_{j}(t) \cong \displaystyle \frac{1}{\Delta x}\int_{x_{j+\frac{1}{2}}}^{x_{j-\frac{1}{2}}} \bu(x,t) \,dx$ is the approximation to the cell average of the solution at the $j$-th cell in time $t$;
	\item $P^t_{j}(x) = P_j(x; \{ \bu_i(t) \}_{i \in  \mathcal{S}_j})$ is a high-order reconstruction operator, i.e. an operator that gives a smooth high-order approximation of the solution at the $j$-th cell from the values of the cell-average approximations available at cells belonging to the stencil $\mathcal{S}_j$;
	\item $\mathcal{D}^\pm_{j+1/2}(t) =\mathcal{D}^\pm_{j+1/2}(\bu^-_{j+1/2}(t),\bu^+_{j+1/2}(t))$, where:
	$$
	\bu_{j+\frac{1}{2}}^{-} (t)= P^t_{j}(x_{j+\frac{1}{2}}),\ \ \bu_{j+\frac{1}{2}}^{+}(t) = P^t_{j+1}(x_{j+\frac{1}{2}})
	$$
	are the reconstructions at the cell interface
	and $\mathcal{D}^\pm(\bu_{l}, \bu_{r})$ are the fluctuations corresponding to a first-order path-conservative numerical methods \eqref{scheme}: see \cite{castro2017well} for details.
 \end{itemize}

The usual strategy consisting in applying an ODE solver to discretize \eqref{eq:hopc}  in time makes difficult the combination of the high-order method with the in-cell discontinuous reconstruction technique that will be described in next section. Therefore, we consider high-order methods based on the integration of (\ref{eq:hopc}) in the interval $[t^{n}, t^{n+1}]$: 
$$ \bu_{j}^{n+1} = \bu_j^n -\frac{1}{\Delta x}\int_{t^{n}}^{t^{n+1}}\left(\mathcal{D}_{j+\frac{1}{2}}^{-}(t) + \mathcal{D}_{j-\frac{1}{2}}^{+}(t) + \int_{x_{j-\frac{1}{2}}}^{x_{j+\frac{1}{2}}}\mathcal{A}(P_{j}^{n}(x,t)) \partial_x P_{j}^{n}(x,t)dx\right)dt,$$
in which the  Cauchy-Kowaleskaya strategy will be followed to compute the reconstruction operators, similar to ADER approach (see \cite{titarev2005ader}, \cite{dumbser2009ader}).

In practice, two quadrature formulas
$$
\int_{x_{j-1/2}}^{x_{j+1/2}} f(x) \, dx \approx \Delta x \sum_{m = 1}^{M_{x}} \beta_{m}^{x} f(x_j^{m}), \quad
\int_{t^{n}}^{t^{n+1}} g(t) \, dt \approx \Delta t \sum_{l = 1}^{M_{t}} \beta_{l}^{t} g(t^{n,l})
$$
are used to approximate the integrals in space and time, respectively, with the desired order of accuracy: 
\begin{eqnarray}\label{fd_pc_method}
    \bu_{j}^{n+1} & = & \bu_{j}^{n} - \frac{\Delta t}{\Delta x}\sum_{l=1}^{M_{t}} \beta_{l}^{t}\left(\mathcal{D}_{j+\frac{1}{2}}^{-}(t^{n,l}) + \mathcal{D}_{j-\frac{1}{2}}^{+}(t^{n,l})\right) \\ \nonumber & & - \Delta t \sum_{l=1}^{M_{t}} \sum_{m=1}^{M_{x}} \beta_{l}^{t} \beta_{m}^{x}\mathcal{A}(P_{j}^{n}(x_{j}^{m},t^{n,l}))
    \partial_x P_{j}^{n}(x_{j}^{m},t^{n,l}).
\end{eqnarray}



\subsubsection{Second-order method}
Using the equations, the first-degree Taylor polynomial of the exact solution can be written as follows:
\begin{eqnarray*}
\mathcal{P}^n_j(x,t) & = & \bu(x_j,t^n) +  \partial_x\bu(x_j, t^n) (x-x_{j})+ \partial_t\bu(x_j, t^n) (t-t^{n}) \\
& = & \bu(x_j,t^n) +  \partial_x\bu(x_j, t^n) (x-x_{j})  -\mathcal{A}(\bu(x_j, t^n))\partial_x \bu(x_j, t^n)(t-t^{n}).
\end{eqnarray*}
The selected reconstruction operator is then the following approximation of this polynomial:
$$P_{j}^{n}(x,t)= \bu_{j}^{n} + \widetilde{\partial_x \bu}_{j}^{n}(x-x_{j})- \mathcal{A}(\bu_{j}^{n})\widetilde{\partial_x \bu}_{j}^{n}(t-t^{n}),$$
where 
$$
\widetilde{\partial_x\bu}_{j}^{n} = \frac{\bu_{j+1}^{n}- \bu_{j-1}^{n}}{2\Delta x}.
$$

Note that as we are going to use MOOD strategy, no limiting is considered here, and we can use centered approximations. 


The mid-point rule is selected both in time and space, so that  \eqref{fd_pc_method}  reduces to 
\begin{eqnarray*}
    \bu_{j}^{n+1} &=& \bu_{j}^{n} - \frac{\Delta t}{\Delta x}\left(\mathcal{D}_{j+\frac{1}{2}}^{-}(t^{n+\frac{1}{2}}) + \mathcal{D}_{j-\frac{1}{2}}^{+}(t^{n+\frac{1}{2}})\right) - \Delta t \left( \mathcal{A}(P_{j}^{n}(x_{j},t^{n+\frac{1}{2}}))
    \partial_x P_{j}^{n}(x_{j},t^{n+\frac{1}{2}})\right),
\end{eqnarray*}
where 
$$\partial_x P_{j}^{n}(x,t) = \widetilde{\partial_x\bu}_{j}^{n}$$
and 
$$
t^{n+\frac{1}{2}} = t^n + \frac{\Delta t}{2}.
$$

\subsubsection{Third-order method}
 The reconstruction operator is, in this case, an approximation of the second-degree Taylor polynomial of the exact solution of the form: 
        $$P_{j}^{n}(x,t)= \widetilde{\bu}_{j}^{n} + \widetilde{\partial_x\bu}_{j}^{n}(x-x_{j})+ \widetilde{\partial_t \bu}_{j}^{n}(t-t^{n}) + \widetilde{\partial^2_x\bu}_{j}^{n}\frac{(x-x_{j})^2}{2} + \widetilde{\partial_x \partial_t \bu}_{j}^{n}(t-t^{n})(x-x_{j}) + \widetilde{\partial^2_t\bu}_{j}^{n}\frac{(t-t^{n})^2}{2},$$
where, in order to achieve third-order accuracy, the following orders of approximation are necessary:
\begin{eqnarray*}
 \widetilde{\bu}_{j}^{n} & \approx &  \bu(x_{j}, t^{n}) + O(\Delta x^3), \\
 \widetilde{\partial_x\bu}_{j}^{n}  & \approx &  \partial_x\bu(x_j, t^n) + O(\Delta x^2),
 \\
 \widetilde{\partial_t\bu}_{j}^{n}  & \approx &  \partial_t\bu(x_j, t^n) + O(\Delta t^2),
 \\
 \widetilde{\partial^2_x\bu}_{j}^{n}  & \approx &  \partial^2_x\bu(x_j, t^n) + O(\Delta x),\\
 \widetilde{\partial^2_t\bu}_{j}^{n}  & \approx &  \partial^2_t\bu(x_j, t^n) + O(\Delta t),\\
 \widetilde{\partial_x\partial_t\bu}_{j}^{n}  & \approx &  \partial_t\partial_x\bu(x_j, t^n) + O(\Delta x + \Delta t).
 \end{eqnarray*}
These approximations are computed as follows:
$$
\widetilde{\bu}_{j}^{n} =  Q^3_j(x_{j}),
$$
where $Q^3_j$ is the interpolation polynomial satisfying
$$
\frac{1}{\Delta x}\int_{x_{j+l-\frac{1}{2}}}^{x_{j+l+\frac{1}{2}}} Q^3_j(x) \, dx = \bu_{j+l}^n, \quad l = -1,0,1,
$$
and
\begin{eqnarray*}
\widetilde{\partial_x\bu}_{j}^{n} & = &  \frac{\widetilde{\bu}_{j+1}^{n}- \widetilde{\bu}_{j-1}^{n}}{2\Delta x}; \\
\widetilde{\partial^2_x\bu}_{j}^{n}  &= & \frac{\widetilde{\bu}_{j+1}^{n}-2\widetilde{\bu}_{j}^{n} + \widetilde{\bu}_{j-1}^{n}}{\Delta x^2}.
\end{eqnarray*}

To approximate the remaining derivatives, the system is first differentiated  with respect to $x$:
\begin{equation}
    \partial_{x}\partial_{t}\bu =  -\partial_x \left( \mathcal{A}(\bu)\right) \partial_x \bu - \mathcal{A}(\bu) \partial^2_x \bu, 
\end{equation}
and with respect to $t$:
\begin{equation}
    \partial^2_{t}\bu = -\partial_t \left( \mathcal{A}(\bu)\right) \partial_x \bu + \mathcal{A}(\bu)\left(\partial_x \left( \mathcal{A}(\bu)\right) \partial_x \bu\right) + \mathcal{A}(\bu)^2 \partial^2_x \bu.
\end{equation}
The derivatives of $\mathcal{A}(\bu)$ are developed as follows:
\begin{itemize}
    \item $\displaystyle \partial_x \left( \mathcal{A}(\bu)\right) =\sum_{k=1}^{N}  \partial_{x} u_{k} \partial_{u_{k}} \mathcal{A}(\bu)$, where
    $$ \partial_{u_{k}} \mathcal{A}(\bu) = \left(\begin{array}{ccc}
        \partial_{u_{k}} a_{1,1}(\bu) &  \cdots & \partial_{u_{k}} a_{1,N}(\bu) \\
         \vdots & \ddots & \vdots \\
         \partial_{u_{k}}a_{N,1}(\bu) &  \cdots & \partial_{u_{k}} a_{N,N}(\bu)
    \end{array}\right)$$.
    \item $\displaystyle \partial_t \left( \mathcal{A}(\bu)\right) = \sum_{k=1}^{N} \partial_{t} u_{k} \partial_{u_{k}} \mathcal{A}(\bu)
     = -\sum_{k,l=1}^{N} a_{k,l}(\bu)  \partial_x {u_{l}} \partial_{u_{k}} \mathcal{A}(\bu).$
    
\end{itemize}
Using these equalities and the above approximations, the final expression of the reconstruction operator $P_{j}^{n}(x,t)$ is as follows:
\begin{eqnarray*}
P_{j}^{n}(x,t) & = & \widetilde{\bu}_{j}^{n} + \widetilde{\partial_x\bu}_{j}^{n}(x-x_{j})- \mathcal{A}(\widetilde{\bu}_{j}^{n})\widetilde{\partial_x\bu}_{j}^{n}(t-t^{n}) + \widetilde{\partial^2_x\bu}_{j}^{n}\frac{(x-x_{j})^2}{2} \\ & &  -\left(\sum_{k=1}^{N} \left(\widetilde{\partial_{x} \bu}_{j,k}^{n} \partial_{\bu_{j,k}^{n}} \mathcal{A}(\widetilde{\bu}_{j}^{n})\right) \widetilde{\partial_x \bu}_{j}^{n} + \mathcal{A}(\widetilde{\bu}_{j}^{n}) \widetilde{\partial^2_x \bu}_{j}^{n}\right)(t-t^{n})(x-x_{j}) \\ & & +\left(\mathcal{A}(\widetilde{\bu}_{j}^{n})\sum_{k=1}^{N} \left(\widetilde{\partial_{x} \bu}_{j,k}^{n} \partial_{\bu_{j,k}^{n}} \mathcal{A}(\widetilde{\bu}_{j}^{n})\right) \widetilde{\partial_x \bu}_{j}^{n} + \mathcal{A}(\widetilde{\bu}_{j}^{n})^2 \widetilde{\partial^2_x \bu}_{j}^{n} \right. \\
& & \qquad \qquad \left. -\left( \sum_{k,l=1}^{N} a_{k,l}(\widetilde{\bu}_{j}^{n})  
\widetilde{\partial_x \bu}_{j,l}^n \partial_{u_{k}} \mathcal{A}(\widetilde{\bu}_{j}^{n}) \right) \widetilde{\partial_x \bu}_{j}^{n} \right)\frac{(t-t^{n})^2}{2}.
\end{eqnarray*}

The two-point Gauss quadrature formula is used both in time and space, so that  \eqref{fd_pc_method}  writes now as follows: 
\begin{eqnarray*}
    \bu_{j}^{n+1} = \bu_{j}^{n} &-& \frac{1}{2}\frac{\Delta t}{\Delta x}\left(\mathcal{D}_{j+\frac{1}{2}}^{-}(t^{n,1}) + \mathcal{D}_{j-\frac{1}{2}}^{+}(t^{n,1}) + \mathcal{D}_{j+\frac{1}{2}}^{-}(t^{n,2}) + \mathcal{D}_{j-\frac{1}{2}}^{+}(t^{n,2})\right) \\ &+ &  \frac{\Delta t}{4} \Big(\mathcal{A}(P_{j}^{n}(x_{j,1},t^{n,1}))
    \partial_x P_{j}^{n}(x_{j,1},t^{n,1}) + \mathcal{A}(P_{j}^{n}(x_{j,1},t^{n,2}))
    \partial_x P_{j}^{n}(x_{j,1},t^{n,2}) \\ & & \qquad  \mathcal{A}(P_{j}^{n}(x_{j,2},t^{n,1}))
    \partial_x P_{j}^{n}(x_{j,2},t^{n,1}) + \mathcal{A}(P_{j}^{n}(x_{j,2},t^{n,2}))
    \partial_x P_{j}^{n}(x_{j,2},t^{n,2})\Big),
\end{eqnarray*}
where 
$$
x_{j,1} = x_{j-\frac{1}{2}} + \frac{\Delta x}{2}\left(-\frac{1}{\sqrt{3}}+1\right),\ \
x_{j,2} = x_{j-\frac{1}{2}} + \frac{\Delta x}{2}\left(\frac{1}{\sqrt{3}}+1\right).
$$
$$
t^{n,1} = t^{n} + \frac{\Delta t}{2}\left(-\frac{1}{\sqrt{3}}+1\right),\ \
t^{n,2} = t^{n} + \frac{\Delta t}{2}\left(\frac{1}{\sqrt{3}}+1\right).
$$

Following \cite{castro2006} it can be proved that the scheme gives a third-order approximation.


\subsection{MOOD procedure}\label{subsec:mood_procedure}
An unlimited high-order method and a standard first-order path-conservative scheme, whose fluctuations will be represented by
$\mathcal{D}_{j+\frac{1}{2}}^{\pm,LO}$, have to be selected. Observe that another first-order path-conservative scheme has to be selected to compute the 
intercell contributions in the high-order method which is not necessarily the same: its fluctuations will be represented by $\mathcal{D}_{j+\frac{1}{2}}^{\pm,HO}$. The fluctuations corresponding to the Rusanov method will be considered here for $\mathcal{D}_{j+\frac{1}{2}}^{\pm,HO}$.

Once the solutions at time $t^n$ have been computed, the numerical solutions are updated giving the following steps:
\begin{enumerate}
    \item \textbf{Prediction step}: compute
    \begin{equation}\label{CFL}
    \Delta t^{n,CFL} = CFL \min\left( \frac{ \Delta x}{ \max_{j,l}|\lambda^n_{j,l}|} \right)
\end{equation}
where $CFL \in (0,1)$ is the stability parameter and $\lambda^n_{j,l}, \dots, \lambda^n_{j,N}$ represent  the eigenvalues of
$\mathcal{A}(\bu_j^n)$. Advance in time to $t^{n+1} = t^{n}+ \Delta t^{n,CFL}$ using the unlimited high-order method to obtain the candidate solution $\bar{\bu}_{j}^{n+1}$. 

\item \textbf{Cell marking:} following \cite{clain2011high} we consider the relaxed Discrete Maximum Principle (DMP) as our detector: a cell $j$ is marked if 
\begin{equation}\label{DMP}
    \min_{i\in\mathcal{S}_{j}}\{\bu_{i}^{n}\}-\delta_{j} \leq \bar{\bu}_{j}^{n+1} \leq \max_{i\in\mathcal{S}_{j}}\{\bu_{i}^{n}\} + \delta_{j},
\end{equation}
is not satisfied, where 
$$\delta_{j} = \max\left\{tol_1, tol_2\left(\max_{i\in\mathcal{S}_{j}}\{\bu_{i}^{n}\}-\min_{i\in\mathcal{S}_{j}}\{\bu_{i}^{n}\}\right)\right\}.$$
Here, we consider $tol_1 = 10^{-8} $, $tol_2 = 10^{-7}$.
Therefore we define the set of marked cells:
\begin{equation}
    \mathcal{O}_{n} = \{\text{$j$ s.t. \eqref{DMP} is not satisfied}\}.
\end{equation}
More sophisticate marking procedures can also be applied, see for instance \cite{diot2012improved}. The boundary cells, i.e, the cells $j$ such that $j\notin \mathcal{O}_{n}$ but $j-1 \in \mathcal{O}_{n}$ or $j+1 \in \mathcal{O}_{n}$, obtaining thus the set of boundary indices are also marked:
$$\mathcal{B}_{n} = \{\text{ $j \notin \mathcal{O}_{n}$ s.t. $j-1 \in \mathcal{O}_{n}$ or $j+1 \in \mathcal{O}_{n}$} \}$$
\begin{rem}
Observe that $\mathcal{B}_{n} \cap \mathcal{O}_{n} = \emptyset$.
\end{rem}
\item \textbf{Correction step}:
\begin{enumerate}
    \item If $j\notin \mathcal{O}_{n}$ and $j \notin \mathcal{B}_n$, the predicted solution is kept, i.e., $\bu_{j}^{n+1}=\bar{\bu}_{j}^{n+1}$.
    \item If $j \in \mathcal{O}_n$, the selected first-order method is used to compute $\bu_{j}^{n+1}$:
    $$
\bu_j^{n+1}=\bu_j^n-\frac{\Delta t}{\Delta x}\big(\mathcal{D}_{j-1/2}^{+, LO} +
\mathcal{D}_{j+1/2}^{-,LO}\big).
$$
\item If $j \in \mathcal{B}_n$, $\bar \bu_j^{n+1}$ is corrected so that the numerical method is formally consistent with the selected family of paths.
\end{enumerate}
\end{enumerate}

Let us describe more in detail step 3(c): consider a cell $j\in \mathcal{B}_{n}$ such that $j+1\in \mathcal{O}_{n}$. 
In the computation of $\bu_{j+1}^{n + 1}$ using the first-order method,  the right contribution of the nonconservative product at the intercell $x_{j + 1/2}$ is:
$$
\Delta t \, \mathcal{D}_{j+\frac{1}{2}}^{+, LO} = \Delta t \, \mathcal{D}^{+, LO} \left(\bu^-_{j+1/2}(t^{n}),\bu_{j+1}^n \right),
$$ 
while in the computation of the candidate solution at the $j$th cell, the left contribution at the same intercell is: 
$$\Delta t \sum_{l=1}^{M_{t}} \beta_{l}^{t} \mathcal{D}^{-, HO}\left(\bu^-_{j+1/2}(t^{n,l}),\bu^+_{j+1/2}(t^{n,l})\right).$$
The idea is therefore to correct this contribution so that the sum of the left and right ones is equal to
$$
\int_0^1\mathcal{A}(\Phi(s;\bu^-_{j+1/2}(t^{n}),\bu_{j+1}^n ))\frac{\partial\Phi}{\partial
s}(s;\bu^-_{j+1/2}(t^{n}),\bu_{j+1}^n )\,ds.
$$
 A first possibility would be obviously to take $\Delta t\,\mathcal{D}_{j+\frac{1}{2}}^{-,LO}$ as the left contribution what would give:
\begin{eqnarray*}
    \bu_{j}^{n+1}  =  \bu_{j}^{n} & - & \frac{\Delta t}{\Delta x}\left[\sum_{l=1}^{M_{t}} \beta_{l}^{t}
\mathcal{D}_{j-\frac{1}{2}}^{+,HO}(t^{n,l}) + \mathcal{D}_{j+\frac{1}{2}}^{-,LO} \right. \\
    \nonumber 
    & + & \left. \Delta x \sum_{l=1}^{M_{t}} \sum_{m=1}^{M_{x}} \beta_{l}^{t} \beta_{m}^{x}\mathcal{A}(P_{j}^{n}(x_{j}^{m},t^{n,l}))
    \partial_x P_{j}^{n}(x_{j}^{m},t^{n,l}) \right],
\end{eqnarray*}
that, taking into account that the quadrature formula is exact for constant functions, can be written in the form
\begin{eqnarray*}
    \bu_{j}^{n+1}  =  \bu_{j}^{n} & - & \frac{\Delta t}{\Delta x}\left[\sum_{l=1}^{M_{t}} 
    \beta_{l}^{t} \left(
    \mathcal{D}^{+, HO}\left(\bu^-_{j-1/2}(t^{n,l}),\bu^+_{j-1/2}(t^{n,l})\right)
    + \mathcal{D}^{-, LO}\left(\bu^-_{j+1/2}(t^{n}),\bu_{j+1}^{n} \right)
    \right) \right. \\
    & + & \left. \Delta x \sum_{l=1}^{M_{t}} \sum_{m=1}^{M_{x}} \beta_{l}^{t} \beta_{m}^{x}\mathcal{A}(P_{j}^{n}(x_{j}^{m},t^{n,l}))
    \partial_x P_{j}^{n}(x_{j}^{m},t^{n,l}) \right].
\end{eqnarray*}
If this formula is compared with \eqref{fd_pc_method} one can see that, in addition to changing the first-order method used to compute the left contributions at the intercell $x_{j+1/2}$, the left and right reconstructions at this intercell have been set to 
$\bu^-_{j+1/2}(t^{n})$ and $\bu_j^{n+1}$ respectively. Taking into account that the numerical solution in the $j$th cell is approximated by the smooth function $P^n_j$ it is clear that the contributions of the nonconservative products corresponding to the jumps between
$P^n_j(x_{j+1/2}, t^{n,l}) = \bu^-_{j+1/2}(t^{n,l})$ and $\bu^-_{j+1/2}(t^{n})$ for $l =1, \dots, M_t$ are missing in this formula. According to the selected family of paths, these contributions are given by
$$
\int_0^1\mathcal{A}(\Phi(s;\bu_{j+\frac{1}{2}}^{-}(t^{n,l}),\bu^-_{j+1/2}(t^{n})))\frac{\partial\Phi}{\partial
s}(s;\bu_{j+\frac{1}{2}}^{-}(t^{n,l}),\bu^-_{j+1/2}(t^{n}))\,ds.$$
Therefore, the numerical solution at the $j$th is finally computed as follows:
\begin{eqnarray}\label{mood_correction_LO}
    \bu_{j}^{n+1}  =  \bu_{j}^{n} & - & \frac{\Delta t}{\Delta x}\left[\sum_{l=1}^{M_{t}} \beta_{l}^{t}\mathcal{D}_{j-\frac{1}{2}}^{+,HO}(t^{n,l}) + \mathcal{D}_{j+\frac{1}{2}}^{-,LO} \right. \\
    \nonumber 
    & + & \Delta x \sum_{l=1}^{M_{t}} \sum_{m=1}^{M_{x}} \beta_{l}^{t} \beta_{m}^{x}\mathcal{A}(P_{j}^{n}(x_{j}^{m},t^{n,l}))
    \partial_x P_{j}^{n}(x_{j}^{m},t^{n,l}) \\ \nonumber & + &  \left. \sum_{l=1}^{M_{t}} \beta_{l}^{t}\int_0^1\mathcal{A}(\Phi(s;\bu_{j+\frac{1}{2}}^{-}(t^{n,l}),\bu^-_{j+1/2}(t^{n})))\frac{\partial\Phi}{\partial
s}(s;\bu_{j+\frac{1}{2}}^{-}(t^{n,l}),\bu^-_{j+1/2}(t^{n}))\,ds\right].
\end{eqnarray}

Observe that this correction can be equivalently written as follows:
    \begin{eqnarray}\label{mood_correction_LO_2}
    \bu_{j}^{n+1}  =  \bu_{j}^{n} & - & \frac{\Delta t}{\Delta x}\left[\sum_{l=1}^{M_{t}} \beta_{l}^{t}\mathcal{D}_{j-\frac{1}{2}}^{+,HO}(l) + \mathcal{D}_{j+\frac{1}{2}}^{-,LO} \right. \\
    \nonumber 
    & + & \Delta x \sum_{l=1}^{M_{t}} \sum_{m=1}^{M_{x}} \beta_{l}^{t} \beta_{m}^{x}\mathcal{A}(P_{j}^{n}(x_{j}^{m},t^{n,l}))
    \partial_x P_{j}^{n}(x_{j}^{m},t^{n,l}) \\ \nonumber & + &  \left. \sum_{l=1}^{M_{t}} \beta_{l}^{t}\left(\mathcal{D}^{-}(\bu_{j+\frac{1}{2}}^{-}(t^{n,l}), \bu^-_{j+1/2}(t^{n})) + \mathcal{D}^{+}(\bu_{j+\frac{1}{2}}^{-}(t^{n,l}), \bu^-_{j+1/2}(t^{n}))\right)\right].
\end{eqnarray}

A similar expression is found for the correction in a cell $j\in \mathcal{B}_n$ such that $j-1\in \mathcal{O}_{n}$.

An important remark is that these corrections lead to a numerical method that is conservative if the system is conservative, i.e. if 
$\mathcal{A}(\bu)$ is the Jacobian of a flux function $F(\bu)$. In effect, let us represent by $\mathcal{F}^{LO}$ and $\mathcal{F}^{HO}$
the consistent numerical flux defined from $\mathcal{D}^{\pm,LO}$ and $\mathcal{D}^{\pm,LO}$ respectively using \eqref{numflux1} or \eqref{numflux2}. Then, one has:

\begin{eqnarray}\label{proof_conservation}
    \bu_{j}^{n+1} & = & \bu_{j}^{n}  - \frac{\Delta t}{\Delta x}\left[ \sum_{l=1}^{M_{t}} \beta_{l}^{t}\left(F(\bu_{j-\frac{1}{2}}^{+}(t^{n,l})) - \mathcal{F}_{j-\frac{1}{2}}^{HO}(l) \right) + \mathcal{F}_{j+\frac{1}{2}}^{LO} -F(\bu^-_{j+1/2}(t^{n}))  \right. \\ \nonumber & + & \sum_{l=1}^{M_{t}} \beta_{l}^{t} \left(F(\bu_{j+\frac{1}{2}}^{-}(t^{n,l})) - F(\bu_{j-\frac{1}{2}}^{+}(t^{n,l}))\right)   \\ \nonumber & + &  \left. \sum_{l=1}^{M_{t}} \beta_{l}^{t} \left( F(\bu^-_{j+1/2}(t^{n}))- F(\bu_{j+\frac{1}{2}}^{-}(t^{n,l}))\right)\right] \\ \nonumber & = & \bu_{j}^{n}  - \frac{\Delta t}{\Delta x}\left[ \mathcal{F}_{j+\frac{1}{2}}^{LO} - \sum_{l=1}^{M_{t}} \beta_{l}^{t}\mathcal{F}_{j-\frac{1}{2}}^{HO}(t^{n,l})   \right].
\end{eqnarray}
Similar calculations lead to the following expressions for $\bu_{j-1}^n$ and $\bu_{j+1}^n$:
$$
 \bu_{j-1}^{n+1}  = \bu_{j-1}^{n} - \frac{\Delta t}{\Delta x}
 \left[ \sum_{l=1}^{M_t}\beta_t \mathcal{F}_{j-\frac{1}{2}}^{HO}(t^{n,l}) - \sum_{l=1}^{M_{t}} \beta_{l}^{t}\mathcal{F}_{j-\frac{3}{2}}^{HO}(t^{n,l})   \right],
$$
$$
 \bu_{j+1}^{n+1}  =  \bu_{j+1}^{n} - \frac{\Delta t}{\Delta x}\left( \mathcal{F}_{j+\frac{3}{2}}^{LO} - \mathcal{F}_{j+\frac{1}{2}}^{LO}\right),
 $$
 what shows that the numerical fluxes at  $x_{j \pm 1/2}$ coincide. In particular, even if the system is nonconservative, this reasoning shows that the numerical method will be conservative for the conservations laws included in the system.

\section{DR.MOOD methods}\label{sec:h-o_DRMOOD}
MOOD methods introduced in the previous section are not expected in general to capture correctly the discontinuities of the weak solutions to be approximated: close to a discontinuity, the unlimited high-order methods will produce oscillations so that the selected standard first-order path-conservative method will deal with it and, as it has been mentioned, errors in the discontinuity location and/or amplitude are expected. In order to overcome this difficulty we propose to use, instead of a standard first-order method, the in-cell discontinuous reconstruction method introduced in \cite{chalons2019path} and extended to second-order in \cite{pimentel2021cell}. The main difficulty comes from the fact that this method uses time steps that are related not only to stability through the CFL condition but also to the propagation of discontinuous waves through the cells. Therefore, the high-order predictor method and the in-cell discontinuous reconstruction first-order schemes will use in general different time steps, what adds an extra difficulty to the numerical treatment of the cells placed at the boundaries between high-order and first-order regions.  

Besides high-order methods, first-order DR.MOOD methods will be also considered here, i.e. methods that use the MOOD strategy to combine a standard path-conservative first-order scheme with the in-cell discontinuous reconstruction one: on the one hand, such a method is expected to correctly capture isolated discontinuous waves; on the other hand, it is expected to be much more efficient than using the in-cell discontinuous reconstruction technique in the whole computational domain. In effect, small time steps related to the evolution of the discontinuous waves through the cells will be only used close to their location, while the time step given by the CFL condition will be used in the rest of the computational domain.

\subsection{First-order in-cell discontinuous reconstruction path-conservative methods}\label{sec:O1_DR}

The numerical solutions provided by the in-cell discontinuous reconstruction technique are updated using a formula similar to the one of high-order methods \eqref{eq:hopc}, namely:
\begin{equation}\label{eq:DR}
    \bu_{j}^{n+1} = \bu_{j}^{n} - \frac{\Delta t_n}{\Delta x}\left(\mathcal{D}_{j+\frac{1}{2}}^{-} + \mathcal{D}_{j-\frac{1}{2}}^{+} + \dashint_{x_{j-\frac{1}{2}}}^{x_{j+\frac{1}{2}}}\mathcal{A}(P_{j}^{n}(x, t^n))\frac{\partial}{\partial x}P_{j}^{n}(x, t^n)dx\right),
\end{equation}
where 
$$\mathcal{D}^\pm_{j+1/2} =\mathcal{D}^\pm_{j+1/2}(\bu^-_{j+ 1/2},\bu^+_{j+1/2}), $$
with
$$\bu^-_{j+1/2} = P_{j}^{n}(x_{j+\frac{1}{2}}, t^n), \quad \bu^+_{j+1/2}= P_{j+1}^{n}(x_{j+\frac{1}{2}}, t^n).$$
Here, $\mathcal{D}^\pm_{j+1/2}$ are the fluctuations of a standard first-order path-conservative method and the $P^n_j(x,t)$ are reconstruction functions that,  in this case, are not used to increase the accuracy of the method but to remove the numerical viscosity in the cells where a discontinuity is detected. Therefore, once the numerical approximations $\bu_j^n$ at time $t^n = n \Delta t$ have been computed, the first step is to mark the cells $I_j$ such that the solution of the Riemann problem consisting of
\eqref{sys:nonconservative} and the initial conditions
\begin{equation} \label{RPj}
\bu(x,0) = \begin{cases} 
\bu^n_{j-1} & \text{if $x<0$,}\\
\bu^n_{j+1} & \text{if $x >0$,}
\end{cases}
\end{equation}
involves a shock wave or a contact discontinuity.  Let us denote by $\mathcal{M}_n$ the set of indices of the marked cells. 
Then, the piecewise constant reconstructions  $P_{j}^{n}(x,t)$ are defined as follows:
    \begin{itemize}
\item If $j \in \mathcal{M}_n$ then
        $$P_{j}^{n}(x,t) = 
\begin{cases}
\bu^n_{j,l} & \text{ if $x \leq x_{j - 1/2} + d_{j}^{n} \Delta x + \sigma_j^n (t -t^n) $,}\\
\bu^n_{j,r} & \text{otherwise},
\end{cases}.$$
where
 $d_j^n$ is chosen so that
\begin{equation}\label{conservation}
d_j^n u_{j,l, k}^n + (1 -  d_j^n) u_{j,r, k}^n = u_{j, k}^n, 
\end{equation}
for some index $k \in \{ 1, \dots, N \}$; and $\sigma_j^n$, $\bu^n_{j,l}$, and $\bu^n_{j,r}$ are chosen so that if $\bu^n_{j-1}$ and
$\bu^n_{j+1}$ may be  linked by an admissible discontinuity with speed $\sigma$, then
\begin{equation}\label{speedstates}
\bu^n_{j,l} = \bu^n_{j-1}, \quad \bu^n_{j,r} = \bu^n_{j+1}, \quad \sigma_j^n = \sigma.
\end{equation}
See Figure \ref{fig:In_cell_discontinuous_reconstruction_operator}. Observe that this in-cell discontinuous reconstruction can only be done
if $0 \leq d_j^n \leq 1$, i.e. if
$$
0 \leq \frac{u_{j,r,k}^n - u_{j,k}^n}{u_{j,r,k}^n - u_{j,l,k}^n} \leq 1,
$$
otherwise the index $j$ is removed from the set $\mathcal{M}_n$ and constant reconstruction is applied in the cell.
Moreover, if $d_j^n = 1$ and $\sigma_j^n > 0$ (resp. $d_j^n = 0$ and $\sigma_j^n < 0$) the  cell is unmarked and the cell $I_{j+1}$ (resp. $I_{j-1}$) is marked if necessary: note that in these cases, the discontinuity leaves the cell $I_j$ for any $t > t^n$.

\begin{figure}[htpb]
		\begin{subfigure}{0.5\textwidth}
			\includegraphics[width=0.9\linewidth]{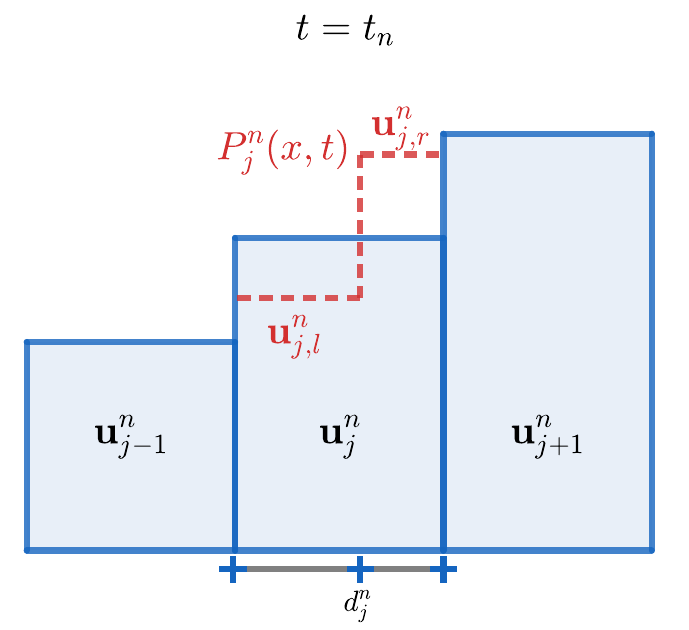}
			\caption*{In-cell discontinuous reconstruction operator at time $t = t^{n}$.}
		\end{subfigure}
		\begin{subfigure}{0.5\textwidth}
				\includegraphics[width=0.9\linewidth]{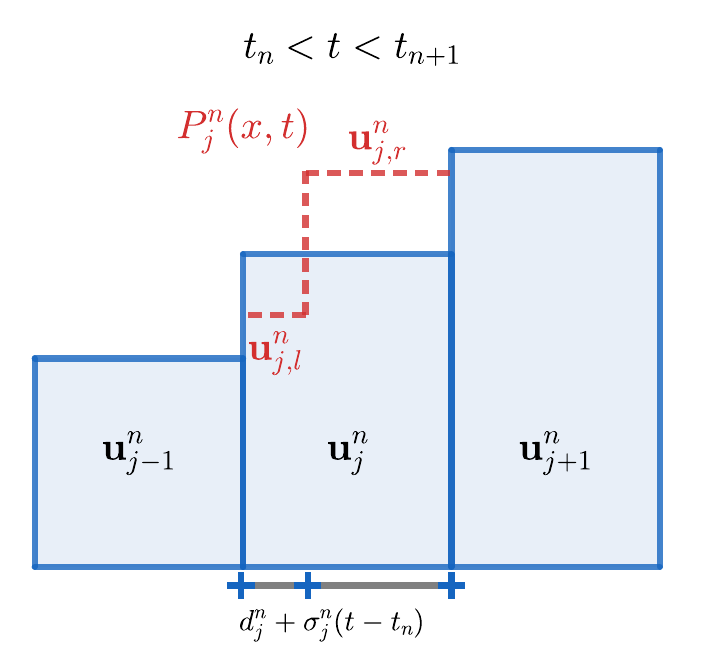}
				\caption*{In-cell discontinuous reconstruction operator at intermediate time $t^{n} < t < t^{n+1}$}
		\end{subfigure}
		\caption{In-cell discontinuous reconstruction operator.}
		\label{fig:In_cell_discontinuous_reconstruction_operator}
	\end{figure}
       
        \item If $j \notin \mathcal{M}_n$ then
        \begin{equation}\label{eq:Constant_reconstruction_operator}
        P_{j}^{n}(x,t)= \bu_{j}^{n}.
        \end{equation}
        See Figure \ref{fig:Constant_reconstruction_operator}.
        \begin{figure}[htpb]
        \centering
			\includegraphics[width=0.5\linewidth]{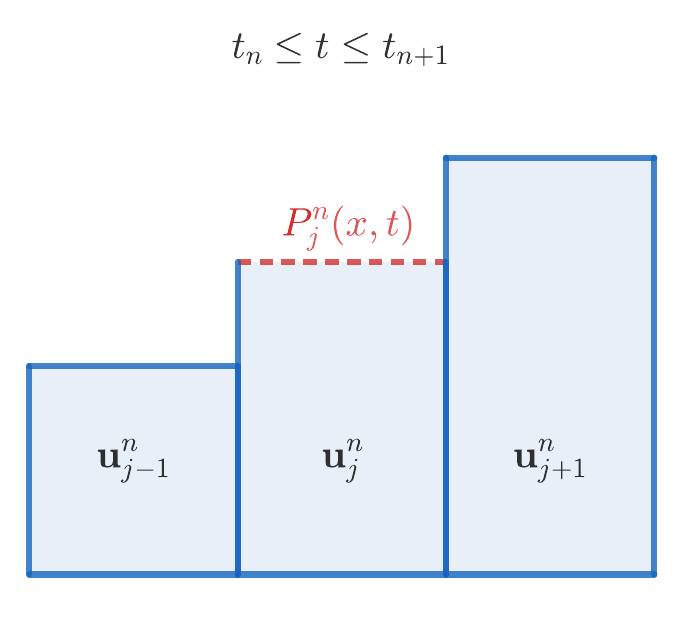}
			\caption{Constant reconstruction operator.}
		\label{fig:Constant_reconstruction_operator}
	\end{figure}
\end{itemize}

\begin{rem} In the case $j \in \mathcal{M}_n$, if one of the equations of system \eqref{sys:nonconservative}, say the $k$th one,  is a conservation law,   the index $k$ is selected in  \eqref{conservation},  so that the corresponding variable is conserved. Moreover, if there is a linear combination of the unknowns $\sum_{k=1}^N \alpha_k u_k$ that is conserved, \eqref{conservation} may be replaced by:
\begin{equation}\label{conservation2}
d_j^n \sum_{k = 1}^N \alpha_k u_{j,l, k}^n + (1 -  d_j^n) \sum_{k = 1}^N \alpha_k u_{j,r, k}^n = \sum_{k = 1}^N \alpha_k u_{j, k}^n. 
\end{equation}
If there are more than one conservation laws, the index $k$ corresponding to one of them is selected in \eqref{conservation}.

\end{rem}

It can be easily checked that the numerical method can be equivalently rewritten as follows:
\begin{equation}\label{eq:discrete1}
    \bu_{j}^{n+1} = \bu_{j}^{n} - \frac{\Delta t_n}{\Delta x}\left(\mathcal{D}_{j+\frac{1}{2}}^{-} + \mathcal{D}_{j-\frac{1}{2}}^{+} + \mathcal{D}_j\right),
\end{equation}
where
\begin{equation}\label{eq:discrete2}
 \mathcal{D}_j = 
 \begin{cases}
\displaystyle  \int_0^1 \mathcal{A}(\Phi(s; \bu^n_{j,l}, \bu^n_{j,r}))\partial_s \Phi(s; \bu^n_{j,l}, \bu^n_{j,r})\, ds & \text{if $j \in \mathcal{M}_n$;}\\
 0 &\text{otherwise.}
 \end{cases}
\end{equation}



\subsubsection{Cell marking} \label{ss:m_n}
Different strategies can be followed to define the set of marked cells $\mathcal{M}_n$ among them: 

\begin{itemize}

\item \textbf{First strategy:}
 If the solutions for Riemann problems are explicitly known, then a cell is marked it the solution of the Riemann problem with initial conditions \eqref{RPj} involves at least a discontinuous wave (contact discontinuity or shock wave). In this case, $\sigma_j^n$, $\bu^n_{j,l}$, $\bu^n_{j,r}$ are selected  as the speed,  the left, and the right states of one of the discontinuous wave. The natural choice for the fluctuations  $\mathcal{D}^{\pm}$ is then Godunov method \eqref{NCFG1}-\eqref{NCFG2}.

\item \textbf{Second strategy:} Let us assume that a Roe linearization $\mathcal{A}_\Phi(\bu_l, \bu_r)$, is available. Then, for every $j$ the coordinates $\{ \alpha_k \}$ of  $\bu^{n}_{j+1} - \bu^{n}_{j-1}$ in the basis of eigenvectors of $\mathcal{A}_\Phi(\bu^{n}_{j-1}, \bu^{n}_{j+1})$, are computed, i.e.
$$
\bu^n_{j+1} - \bu^n_{j-1} = \sum_{k=1}^N \alpha_k R_k(\bu^n_{j-1}, \bu^n_{j+1}),
$$
and the intermediate states
$$
\bu^{(k)}_j = \bu^n_{j-1} + \sum_{s=1}^k \alpha_s R_s (\bu^n_{j-1}, \bu^n_{j+1}),
\quad k = 0, \dots, N,
$$
are defined.
The cell is then marked if there exists an index $k^*$ such that
\begin{equation}\label{eq:choice_k}
    \frac{|\alpha_k^{*}|}{\sum_{k=1}^{N} |\alpha_k|} >\overline{\alpha} \quad \text{and}\quad \lambda_{k^*}(\bu_{j-1}^{n}) \geq \lambda_{k^*}(\bu_{j+1}^{n}),
\end{equation}
where $\overline{\alpha} \in (0,1]$ is a parameter to be selected. Observe that, if \eqref{eq:choice_k} is satisfied with $\overline{\alpha} \approx 1$, then the solution of the linear Riemann problem essentially consists of only one discontinuity that, due to the Roe property is also an admissible discontinuity of the system. Moreover, the second condition in \eqref{eq:choice_k} allows one to discard non-entropy shocks. In this case, the following choices are made:
$$
\sigma_j^n = \lambda_{k^*}(\bu^n_{j-1}, \bu^n_{j+1}), 
\quad \bu^n_{j,l} = \bu^{(k^*-1)}_j, \quad \bu^n_{j,r} = \bu^{(k^*)}_j.
$$ 
Summing up, the set of marked cells is given by
$$
\mathcal{M}_{n} = \{ j \in \mathcal{O}^n \text{ s.t. there exists $k^{*} \in \{1, \dots, N\}$ satisfying \eqref{eq:choice_k}.} \}
$$
This strategy, that is introduced in this paper, will be followed in the numerical tests with $\overline{\alpha}=0.9$. In this case, the natural choice for the fluctuations $\mathcal{D}^{\pm}$  is Roe method \eqref{NCFR}.

\end{itemize}

Observe that, if the solution of the Riemann problem consists of only one discontinuous wave of speed $\sigma$ linking  $\bu^n_{j-1}$ and $\bu^n_{j+1}$, then necessarily \eqref{speedstates} is satisfied for both strategies.





\subsubsection{Time step}\label{ss:timesteps}

The time step $\Delta t_n$ is chosen as follows:
\begin{equation}\label{restricted_time_step}
    \Delta t_n = \min(\Delta t^c_n, \Delta t^r_n).
\end{equation}
Here
\begin{equation}
    \Delta t_n^c = CFL\cdot \min_j\left( \frac{ \Delta x}{ \max_{l}|\lambda_{j,l}|} \right)
\end{equation}
where $CFL \in (0,1)$ is the stability parameter, $\lambda_{j,l}, \dots, \lambda_{j,N}$ represent  the eigenvalues of
$\mathcal{A}(\bu_j^n)$, and
\begin{equation}
    \Delta t_n^r = \min_{j \in \mathcal{M}_n} \begin{cases}
    \displaystyle \frac{1 - d_j^n}{|\sigma_j^n|}\Delta x, & \quad if \quad \sigma_j^n>0, \\ \\
    \displaystyle \frac{d_j^n}{|\sigma_j^n|}\Delta x, & \quad if \quad \sigma_j^n<0.
    \end{cases} 
\end{equation}
Observe that, besides the stability requirement, this choice of time step ensures that no discontinuous reconstruction leaves a marked cell. Note that this time restriction could be avoid as in \cite{chalons2019path} but, for simplicity, we have considered it.

\subsubsection{Shock-capturing property}

We recall the theorem in \cite{chalons2019path, pimentel2021cell} that states that the scheme \eqref{eq:discrete1} captures exactly isolated shocks in the following way:

\begin{thm}\label{th}
Assume that $\bu_l$ and $\bu_r$ can be joined by an entropy shock of speed $\sigma$. Then, the numerical method  \eqref{eq:discrete1} provides an exact numerical solution
of the Riemann problem with initial conditions
$$
\bu(x,0) = \begin{cases}
\bu_l & \text{if $x < 0$,} \\
\bu_r & \text{otherwise,}
\end{cases} 
$$
in the sense that
\begin{equation}\label{presave}
\bu_j^n = \frac{1}{\Delta x}\int_{x_{j-1/2}}^{x_{j+1/2}} \bu(x,t^n) \, dx, \quad \forall j, n
\end{equation}
where $\bu(x,t)$ is the exact solution.

\end{thm}

\subsection{DR.MOOD procedure}\label{DRMOOD_procedure}
 As in MOOD methods, two first-order path-conservative schemes have to be selected: the one used in the in-cell reconstruction method (Godunov or Roe methods, depending on the strategy selected to mark the cells), whose fluctuations will be represented by
$\mathcal{D}_{j+\frac{1}{2}}^{\pm,LO}$, and the one selected to compute the intercell contributions in the high-order method, whose fluctuations will be represented by $\mathcal{D}_{j+\frac{1}{2}}^{\pm,HO}$. Roe and Rusanov methods respectively will be considered in the numerical tests shown in Section \ref{sec:numerical_tests}. Moreover, two different reconstruction operators are involved:   the polynomial reconstruction operator $P^n_j$ used in the unlimited high-order methods to increase the accuracy and the piecewise constant in-cell reconstruction operator $P^{DR,n}_j$.

\begin{figure}[htpb]
        \centering
			\includegraphics[width=0.8\linewidth]{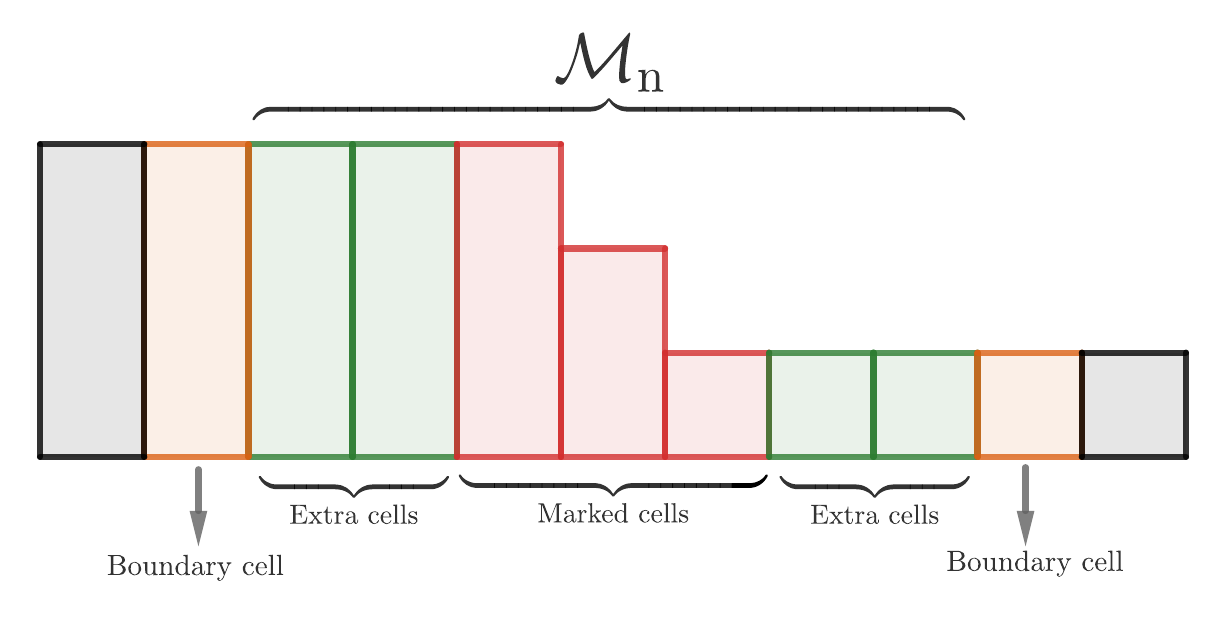}
			\caption{Cell classification.}
		\label{fig:Cell_classification}
	\end{figure}

Once the numerical approximations $\bu_j^n$
of the averages of the solutions have been computed at time $t^n = n \Delta t$, a similar procedure as the one described in Subsection \ref{subsec:mood_procedure} is followed:
\begin{enumerate}
    \item \textbf{Prediction step}: compute $t^{n,CFL}$ \eqref{CFL}  and advance in time to $t^{n+1} = t^{n}+ \Delta t^{n,CFL}$ using the unlimited high-order method to obtain the candidate solution $\bar{\bu}_{j}^{n+1}$.

\item \textbf{Cell marking}: first the cells in which the discrete maximum principle is not satisfied are marked
\begin{equation}\label{oscillating_set}
    \mathcal{O}^{DMP}_{n} = \{\text{$j$ s.t. \eqref{DMP} is not satisfied}\}.
\end{equation}
Then, some of the neighbour of a marked cell are marked as well so that a discontinuity cannot leave the marked region in one time step, namely
   \begin{equation}\label{marked_cells}
    \mathcal{O}_{n} = \mathcal{O}^{DMP}_{n}  \cup \left\{\text{ $j \notin \mathcal{O}^{DMP}_{n}$ s.t. $j\pm s \in\mathcal{O}^{DMP}_{n}$, for any $s \in \left\{1, ..., \frac{|\mathcal{S}| + 1}{2}\right\}$ } \right\},
    \end{equation}
    where $|\mathcal{S}|$ the stencils cardinal. Boundary cells are also marked:
    $$\mathcal{B}_{n} = \{\text{ $j \notin \mathcal{O}_{n}$ s.t. $j-1 \in \mathcal{O}_{n}$ or $j+1 \in \mathcal{O}_{n}$} \}.$$
\item \textbf{Correction step}:
\begin{enumerate}
    \item If $j\notin \mathcal{O}_{n}$ and $j \notin \mathcal{B}_n$, the predicted solution is kept, i.e., $\bu_{j}^{n+1}=\bar{\bu}_{j}^{n+1}$.
    \item If $j \in \mathcal{O}_n$, the in-cell discontinuous reconstruction method is used to compute $\bu_{j}^{n+1}$.
\item If $j \in \mathcal{B}_n$, $\bar \bu_j^{n+1}$ is corrected so that the numerical method is formally consistent with the selected family of paths.

\end{enumerate}
\end{enumerate}
Figure \ref{fig:Cell_classification} shows an example of the cell classification in which $|\mathcal{S}| = 3$ so that two extra cells are added in both sides.

\begin{rem}
    In the case of first-order DR.MOOD methods, the (DMP) criterion is not able to properly detect discontinuities. Therefore,  it is replaced by the following Locally Significant Jump (LSJ) criterion (see \cite{biswas2018low}):
    \begin{equation}\label{LSJ_criterion}
    \frac{|\bu_{j+1}^{n}-\bu_{j-1}^{n}|}{2} > \max(|\bu_{j+2}^{n}-\bu_{j+1}^{n}|, |\bu_{j-1}^{n}-\bu_{j-2}^{n}|).
    \end{equation}
\end{rem}

Observe that, in general, several steps of the in-cell discontinuous reconstruction first-order method has to be given in the sub-mesh corresponding to $j \in \mathcal{O}_n$ to go from $t^n$ to $t^{n+1}$. Let us represent by
$$
t^{n}_0 = t^n < t^{n}_1 < \dots < t^{n}_k = t^{n+1}$$
the corresponding intermediate times. The  time steps
$$
\Delta t_{n,i} = t^{n}_{i+1} - t^{n}_i,\quad i = 0, \dots, k-1
$$
are computed following Subsection \ref{ss:timesteps} (the last one is reduced if necessary so that $t^{n,k} = t^{n+1}$).
The first stage at every one of these intermediate steps is to mark the cells in which a discontinuity is selected: this is done by following one of the strategies mentioned in \ref{ss:m_n}: in particular, the second one will be considered here.
As in MOOD methods, the most involved step in the DR.MOOD procedure is 3(c): technical details are given in  \ref{app}.

Since the in-cell discontinuous reconstruction method is used in cells in which a discontinuity is detected, Theorem \ref{th} holds for DR.MOOD methods.

\section{Numerical tests}\label{sec:numerical_tests}

In this section the MOOD and DR.MOOD methods are applied to two nonconservative systems: the Modified Shallow Water system and the Two-layer Shallow Water system. The following numerical methods are compared:

\begin{itemize}
    \item O$k$\_MOOD: $k$-order path-conservative MOOD method using Rusanov fluctuations \eqref{NCFRus}.
    \item O$k$\_DRMOOD: $k$-order path-conservative DR.MOOD method combining the in-cell discontinuous reconstructions with Roe fluctuations \eqref{NCFR} and the MOOD approach with Rusanov fluctuations \eqref{NCFRus}.
    \item O$k$\_DR: $k$-order path-conservative method with in-cell discontinuous reconstructions developed in \cite{pimentel2021cell}.
\end{itemize}

\subsection{Modified Shallow Water system}\label{ss:MSW}

\subsubsection{Equations}

Let us consider the modified Shallow Water system introduced in \cite{Castro2008}: 
\begin{equation} \label{simplifiedSW}
\left\{
\begin{array}{l}
\partial_t h + \partial_x q = 0,\\
\smallskip
\displaystyle\partial_t q + \partial_x \left(\frac{q^2}{h}\right) + qh\partial_x h = 0,
\end{array}
\right.
\end{equation}
where $\bu = [h,q]^t$ belongs to $\Omega = \{\bu \in \mathbb{R}^{2}| \quad 0<q, \ 0<h<(16q)^{1/3}\}$. This system can be written in the form (\ref{sys:nonconservative}) with
$$
\mathcal{A}(\bu) = \left[ \begin{array}{cc} 0  & 1 \\ -u^2+uh^2 & 2u \end{array} \right],
$$
being $u=q/h$. The system is strictly hyperbolic over $\Omega$ with eigenvalues
$$\lambda_{1}(\bu) = u-h\sqrt{u}, \quad \lambda_{2}(\bu) = u+h\sqrt{u},$$
whose characteristic fields, given by the eigenvectors
$$
R_1(\bu) = [1 , u-h\sqrt{u}]^T, \quad R_2(\bu) = [1, u+h\sqrt{u}]^T,
$$
are genuinely nonlinear.

\subsubsection{Simple waves}

Once the family of paths has been chosen, the simple waves of this system are:
\begin{itemize}
    \item 1-rarefaction waves joining states $\bu_l$, $\bu_r$ such that
    $$ h_r < h_l, \quad \sqrt{u_l} +h_l/2 = \sqrt{u_r} +h_r/2, $$
    and 2-rarefaction waves joining states $\bu_l$, $\bu_r$ such that
    $$ h_l < h_r, \quad \sqrt{u_l} -h_l/2 = \sqrt{u_r} -h_r/2. $$
    
    \item 1-shock and 2-shock waves joining states $\bu_l$ and $\bu_r$ such that
    $h_l <  h_r$ or $h_r < h_l$ respectively, that satisfy the jump conditions:
    \begin{eqnarray*}
    \sigma[h] & = & \left[q \right],\\
    \sigma[q] & = & \left[\frac{q^2}{h}\right] + \int_0^1 \phi_q(s; \bu_l, \bu_r) \phi_h(s; \bu_l, \bu_r) 
    \partial_s \phi_h(s; \bu_l, \bu_r) \,ds.
    \end{eqnarray*}
\end{itemize}
If, for instance, the following family of path is chosen:
$$
\phi(s; \bu_l, \bu_r) =  \left[\begin{array}{c}
     \phi_h(s; \bu_l, \bu_r)  \\
     \phi_q(s; \bu_l, \bu_r) 
\end{array}\right] = \left\{\begin{array}{l}
     \left[\begin{array}{c}
     h_l + 2s(h_r-h_l)  \\
     q_l 
\end{array}\right] \quad \text{if $0\leq s \leq \frac{1}{2}$},  \\
\\
     \left[\begin{array}{c}
     h_r  \\
     q_l + (2s-1)(q_r-q_l) 
\end{array}\right] \quad \text{if $\frac{1}{2} \leq s \leq 1$},
\end{array}\right.
$$
the jump conditions reduce to:
    \begin{eqnarray*}
    \sigma[h] & = & \left[q \right],\\
    \sigma[q] & = & \left[\frac{q^2}{h}\right] + q_l  \left[\frac{h^2}{2}\right].
    \end{eqnarray*}
 If this family of paths has been selected, a Roe matrix is given by
$$
\mathcal{A}(\bu_l,\bu_r)= \left[\begin{array}{cc}
   0  &  1\\
   -\bar{u}^{2}+q_l\bar{h}  & 2\bar{u}
\end{array}\right],
$$
where
$$
\bar{u} = \frac{\sqrt{h_l}u_{l} + \sqrt{h_r}u_{r}}{\sqrt{h_l}+\sqrt{h_r}}, \quad \bar{h} = \frac{h_l+h_r}{2}.
$$

\subsubsection{Cell-marking criterion}

In order to use the (DMP) \eqref{DMP} criterion for the second- and third-order schemes or the (LSJ) criterion \eqref{LSJ_criterion} for the first-order scheme, the variable $h$ is selected in the cell marking procedure.

\subsubsection{Numerical tests}

We present different tests to prove the convergence of DR.MOOD methods to the right weak solution and their order of accuracy. We consider a CFL number of 0.5 and the domain $[-1,1]$ in all tests.
        
\subsubsection*{Test 1: Isolated 1-shock}

\begin{figure}[htpb]
		\begin{subfigure}{0.5\textwidth}
			\includegraphics[width=1.1\linewidth]{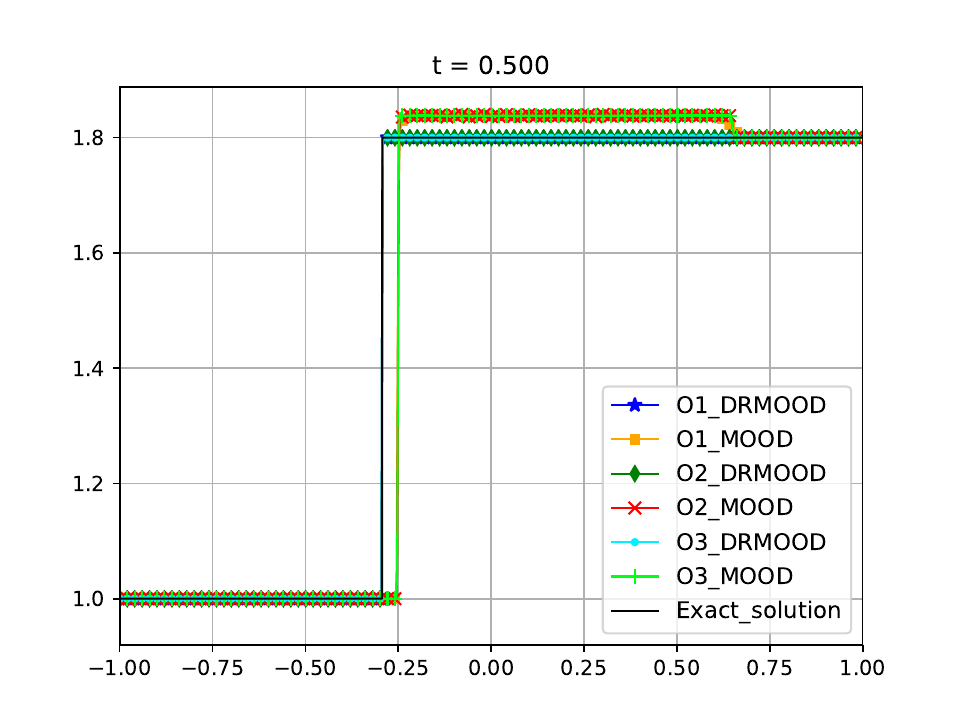}
			\caption{Variable $h$}
		\end{subfigure}
		\begin{subfigure}{0.5\textwidth}
				\includegraphics[width=1.1\linewidth]{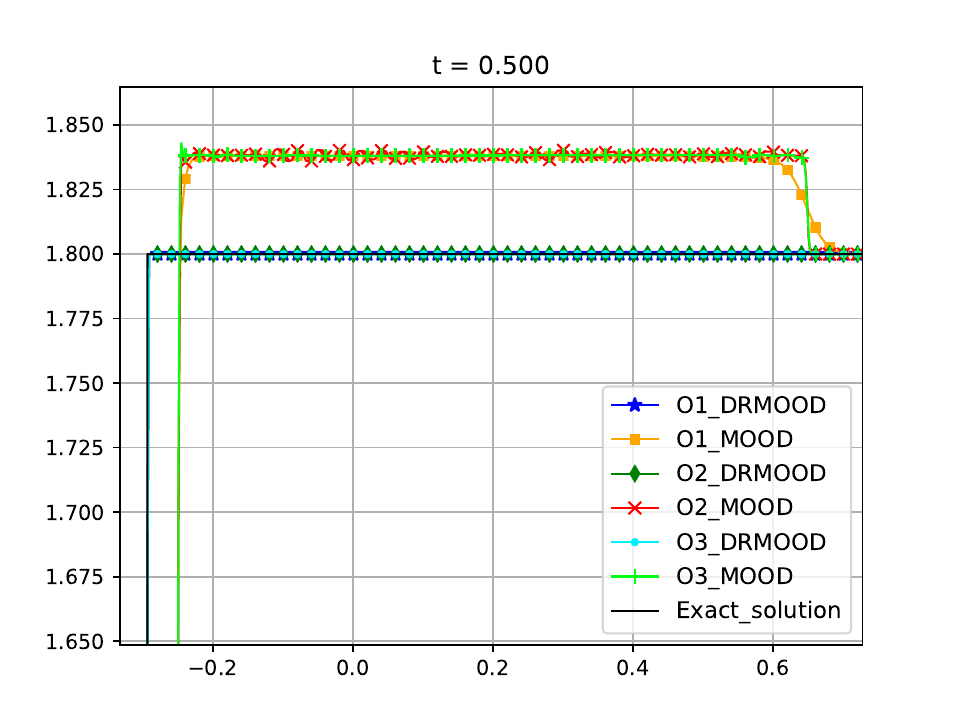}
				\caption*{Variable $h$: Zoom}
		\end{subfigure}
		\begin{subfigure}{0.5\textwidth}
			\includegraphics[width=1.1\linewidth]{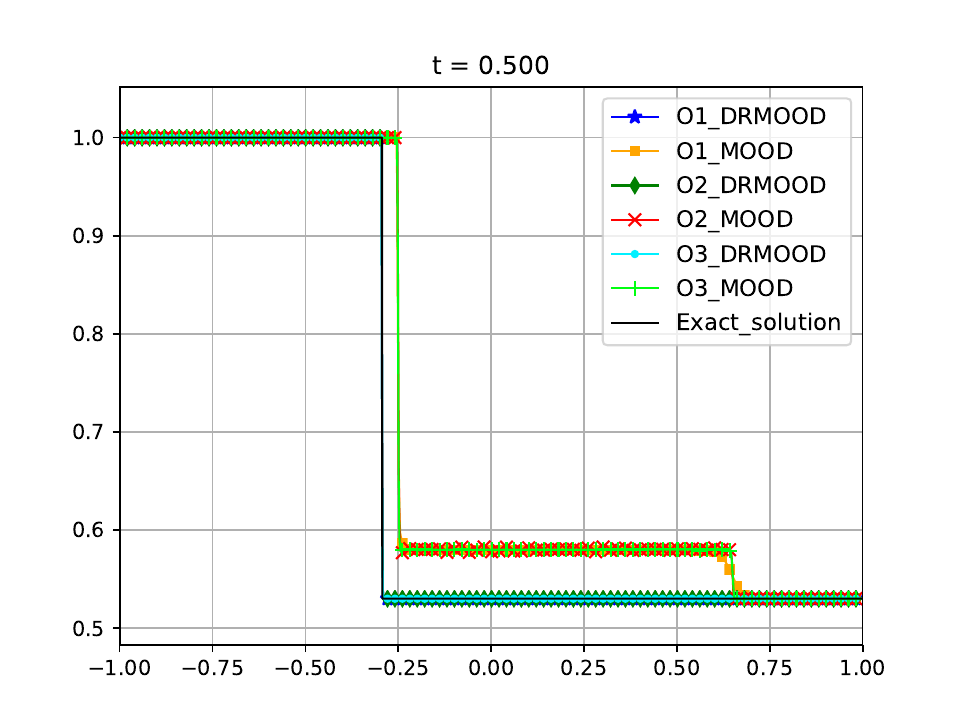}
			\caption{Variable $q$}
		\end{subfigure}
		\begin{subfigure}{0.5\textwidth}
				\includegraphics[width=1.1\linewidth]{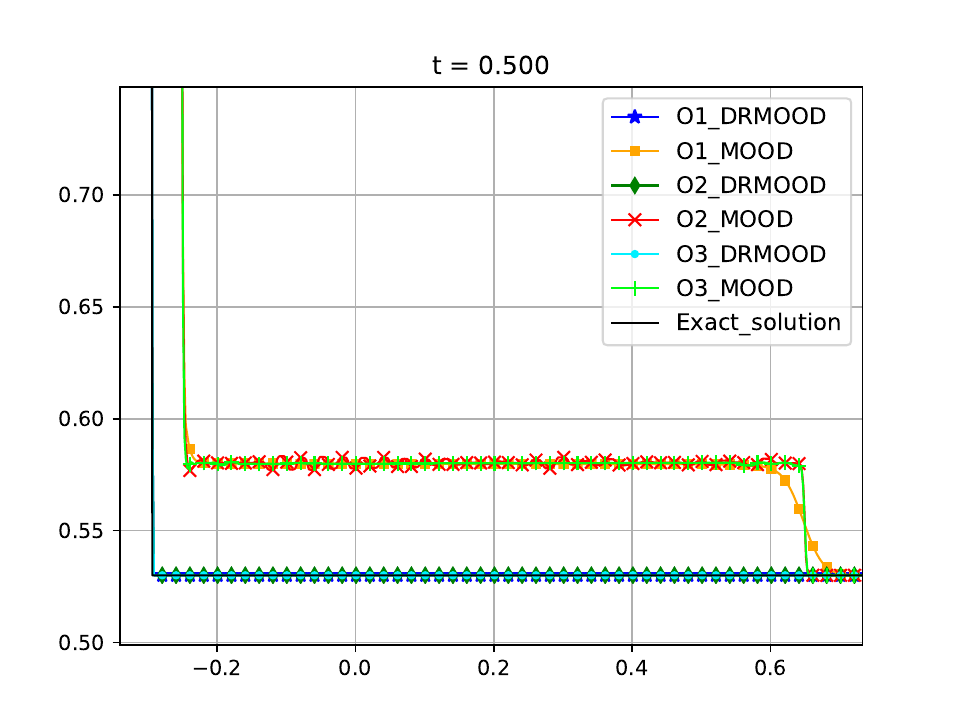}
				\caption*{Variable $q$: Zoom}
		\end{subfigure}
		\caption{Modified Shallow Water system. Test 1: Numerical solutions obtained using the MOOD approach with the first- , second- and third-order methods with and without in-cell discontinuous reconstruction based on the Roe matrix at time $t=0.15$ with 1000 cells. Top: variable $h$ (left), zoom middle state (right). Down: variable $q$ (left), zoom middle state (right).}
		\label{fig:1DModifiedShallowWater_TestIsolated_1000_t05}
	\end{figure}

Let us consider the following initial condition taken from \cite{Castro2008}

$$
\bu_0(x) = [h_0(x),q_{0}(x)]^T= \begin{cases}
     [1, 1]^T & \text{if $x<0$,}  \\
     [1.8,0.530039370688997]^T & \text{otherwise.}
\end{cases}
$$
The solution of the Riemann problem consists of a 1-shock wave joining the left and right states. This test is devoted to show the shock-capturing property of the new methods. In Figure \ref{fig:1DModifiedShallowWater_TestIsolated_1000_t05} we show the numerical results at time $t=0.15$ obtained with O$p$\_MOOD and O$p$\_DRMOOD, $p=1,2,3$, using a 1000-cell mesh. We observe that the only methods that are able to capture well the isolated shock are those with the in-cell discontinuous reconstruction operator. As it was observed in \cite{pimentel2021cell} the error committed by the other methods do not converge to 0 as the mesh is refined. 

\subsubsection*{Test 2: left-moving 1-shock + right-moving 2-shock}\label{ssstest2}

\begin{figure}[htpb]
		\begin{subfigure}{0.5\textwidth}
			\includegraphics[width=1.1\linewidth]{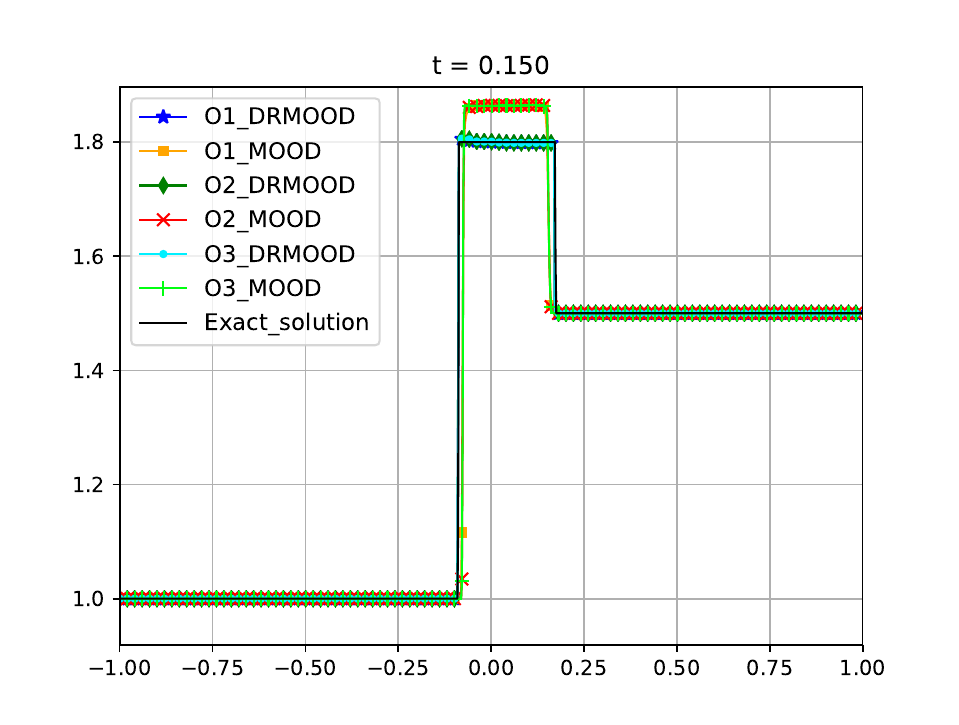}
			\caption{Variable $h$}
		\end{subfigure}
		\begin{subfigure}{0.5\textwidth}
				\includegraphics[width=1.1\linewidth]{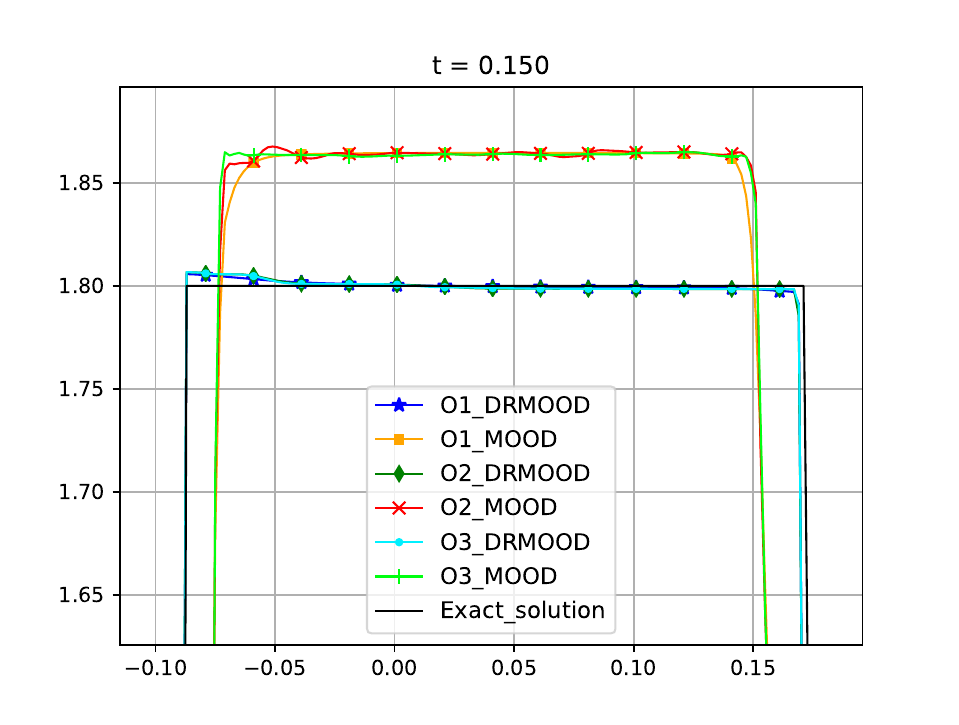}
				\caption*{Variable $h$: Zoom}
		\end{subfigure}
		\begin{subfigure}{0.5\textwidth}
			\includegraphics[width=1.1\linewidth]{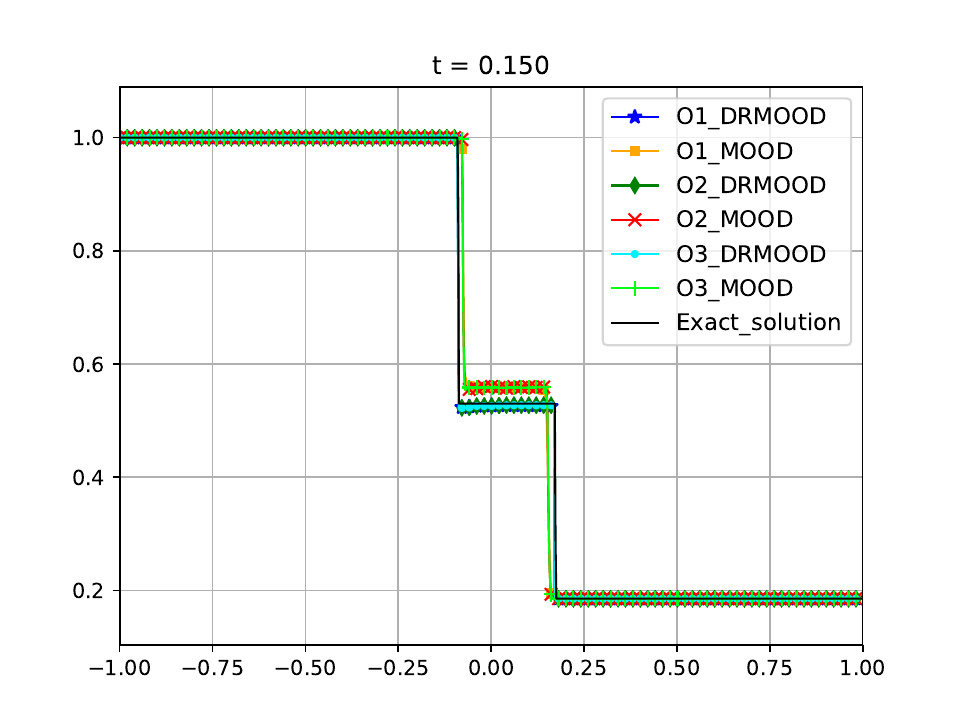}
			\caption{Variable $q$}
		\end{subfigure}
		\begin{subfigure}{0.5\textwidth}
				\includegraphics[width=1.1\linewidth]{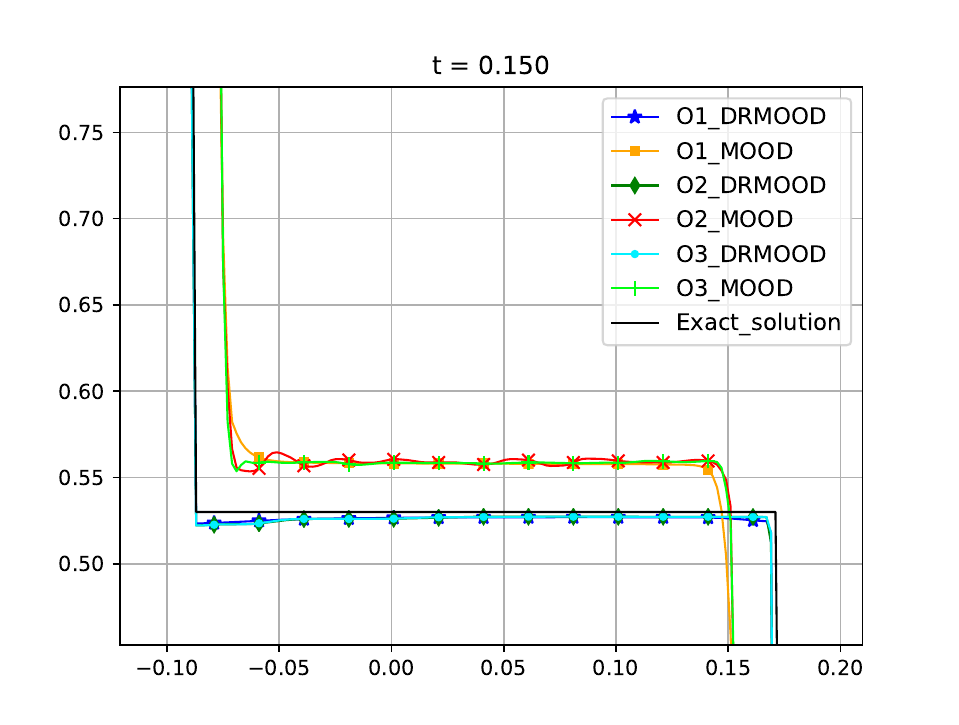}
				\caption*{Variable $q$: Zoom}
		\end{subfigure}
		\caption{Modified Shallow Water system. Test 2: Numerical solutions obtained using the MOOD approach with the first- , second- and third-order methods with and without in-cell discontinuous reconstruction based on the Roe matrix at time $t=0.15$ with 1000 cells. Top: variable $h$ (left), zoom middle state (right). Down: variable $q$ (left), zoom middle state (right).}
		\label{fig:1DModifiedShallowWater_Test1_1000_t015}
	\end{figure}

\begin{figure}[htpb]
		\begin{subfigure}{0.5\textwidth}
			\includegraphics[width=1.1\linewidth]{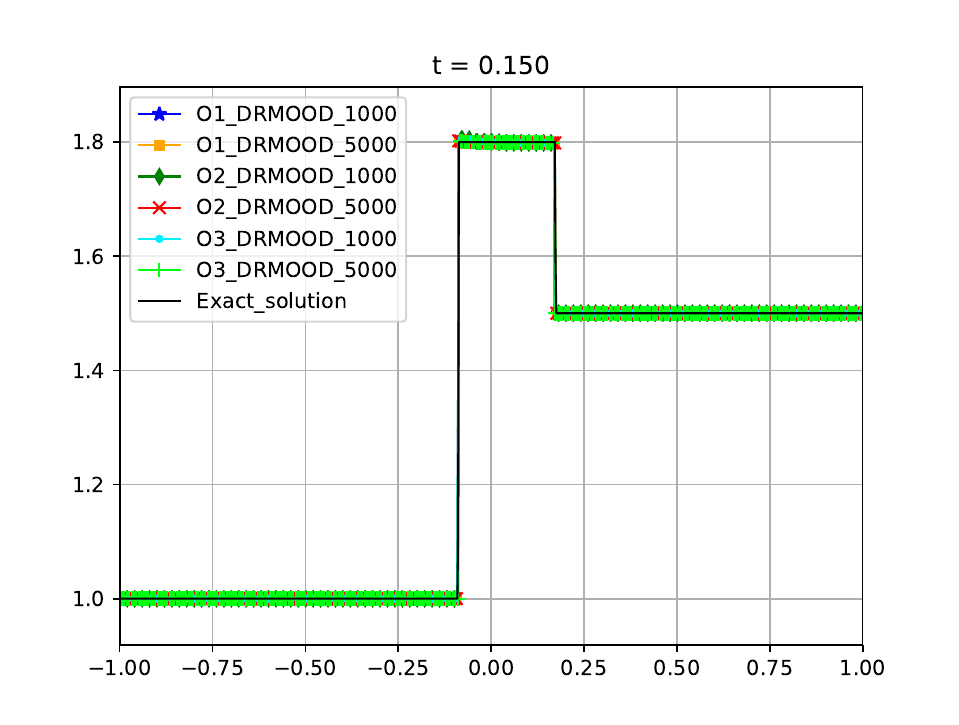}
			\caption{Variable $h$}
		\end{subfigure}
		\begin{subfigure}{0.5\textwidth}
				\includegraphics[width=1.1\linewidth]{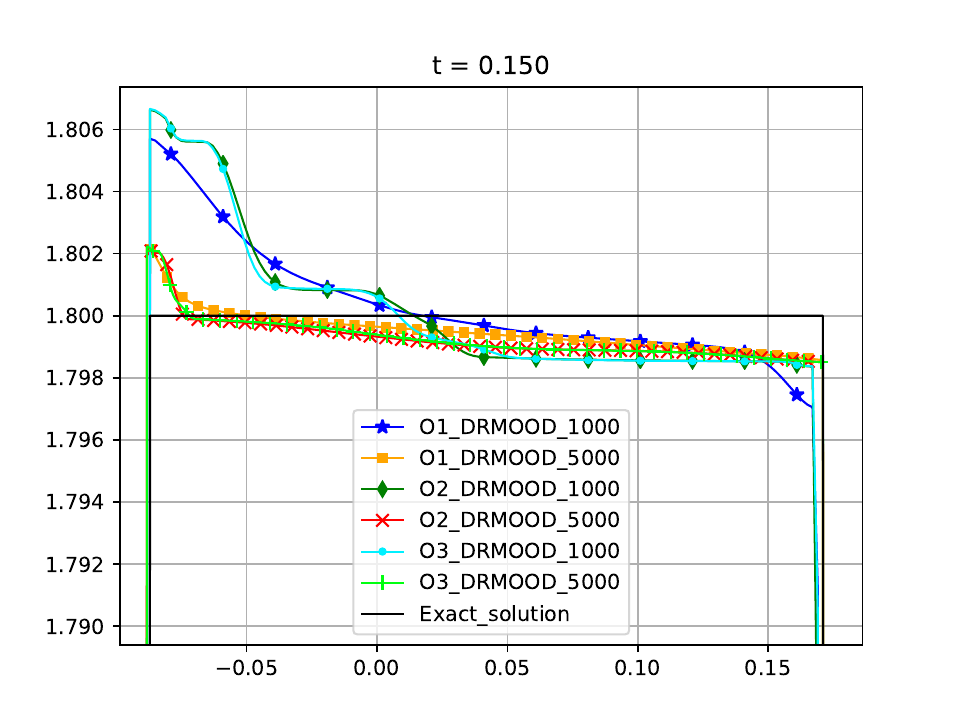}
				\caption*{Variable $h$: Zoom}
		\end{subfigure}
		\begin{subfigure}{0.5\textwidth}
			\includegraphics[width=1.1\linewidth]{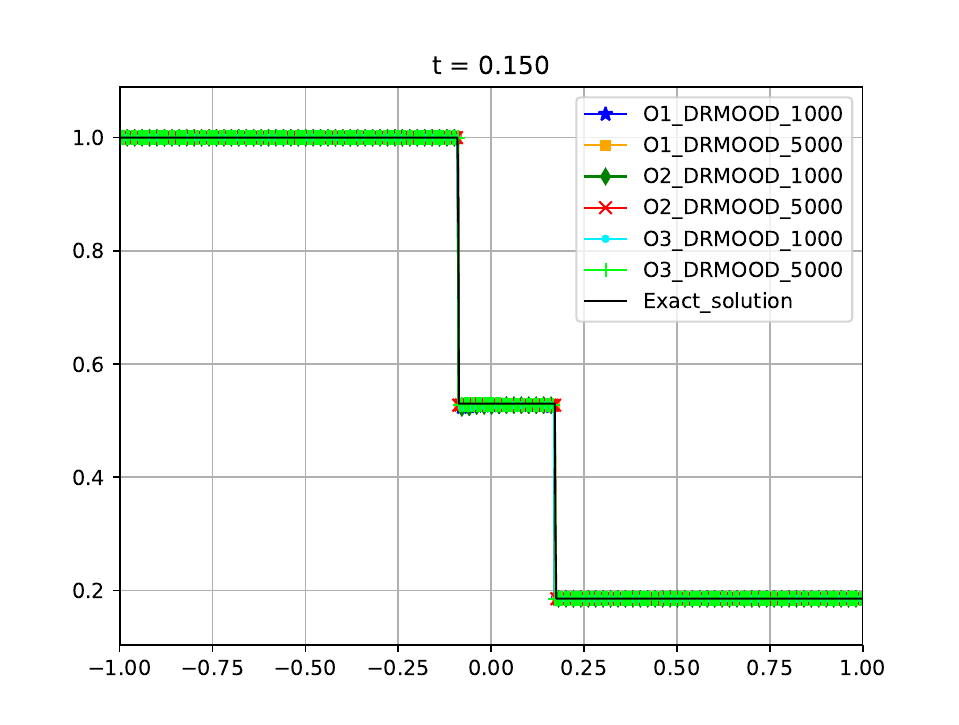}
			\caption{Variable $q$}
		\end{subfigure}
		\begin{subfigure}{0.5\textwidth}
				\includegraphics[width=1.1\linewidth]{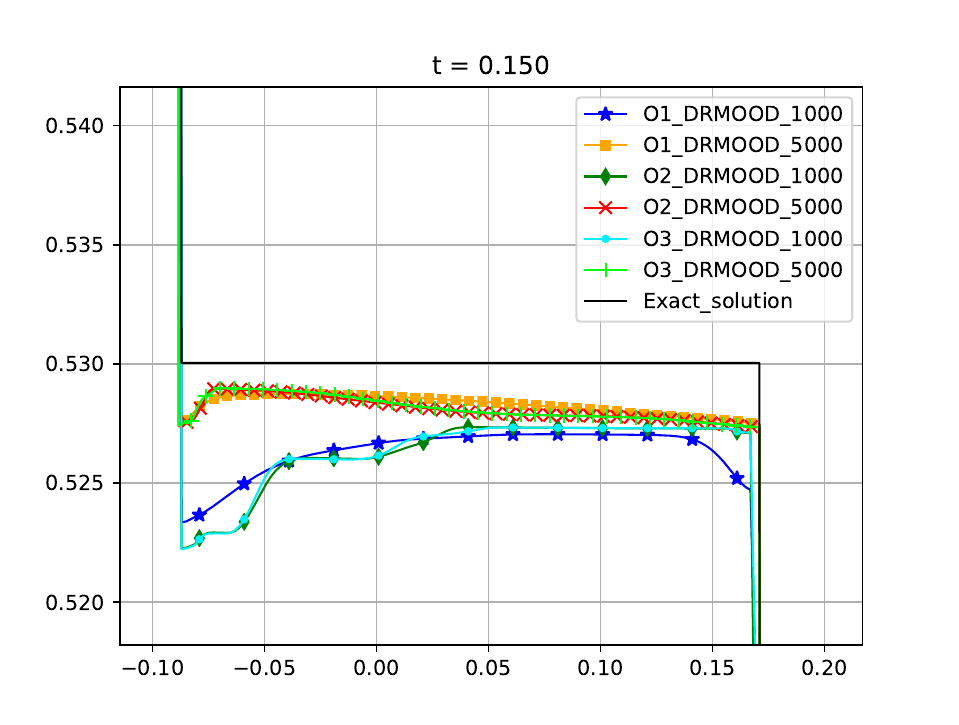}
				\caption*{Variable $q$: Zoom}
		\end{subfigure}
		\caption{Modified Shallow Water system. Test 2: Numerical solutions obtained using 1000- and 5000-cell meshes with the first- , second- and third-order DR.MOOD methods at time $t=0.15$. Top: variable $h$ (left), zoom middle state (right). Down: variable $q$ (left), zoom middle state (right).}
		\label{fig:1DModifiedShallowWater_Test1_1000_vs_5000_t015}
	\end{figure}

\begin{figure}[htpb]
		\begin{subfigure}{0.5\textwidth}
			\includegraphics[width=1.1\linewidth]{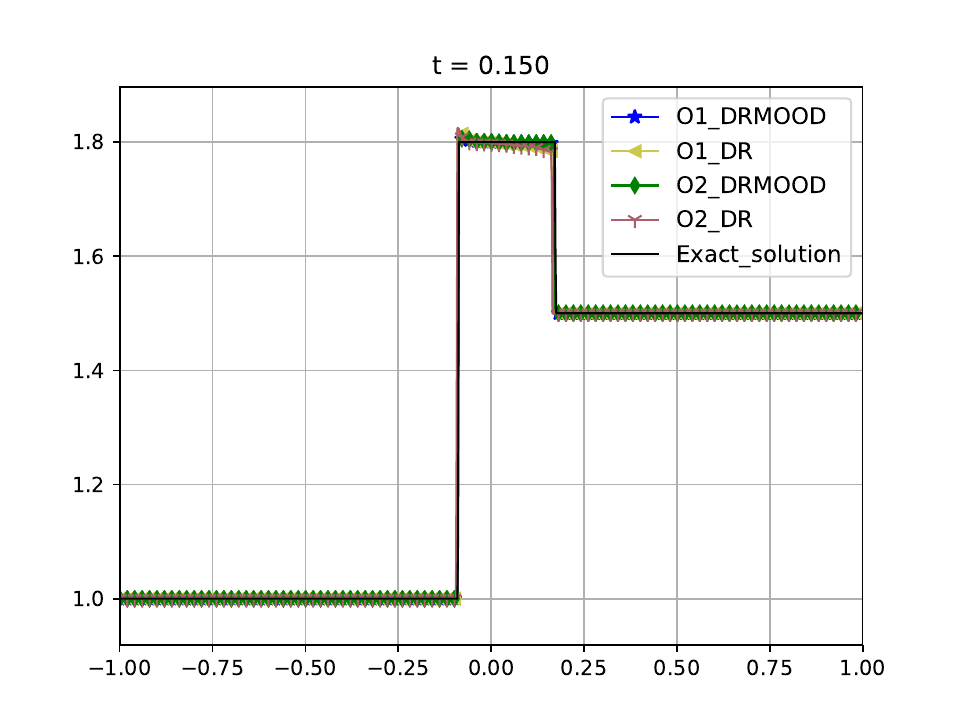}
			\caption{Variable $h$}
		\end{subfigure}
		\begin{subfigure}{0.5\textwidth}
				\includegraphics[width=1.1\linewidth]{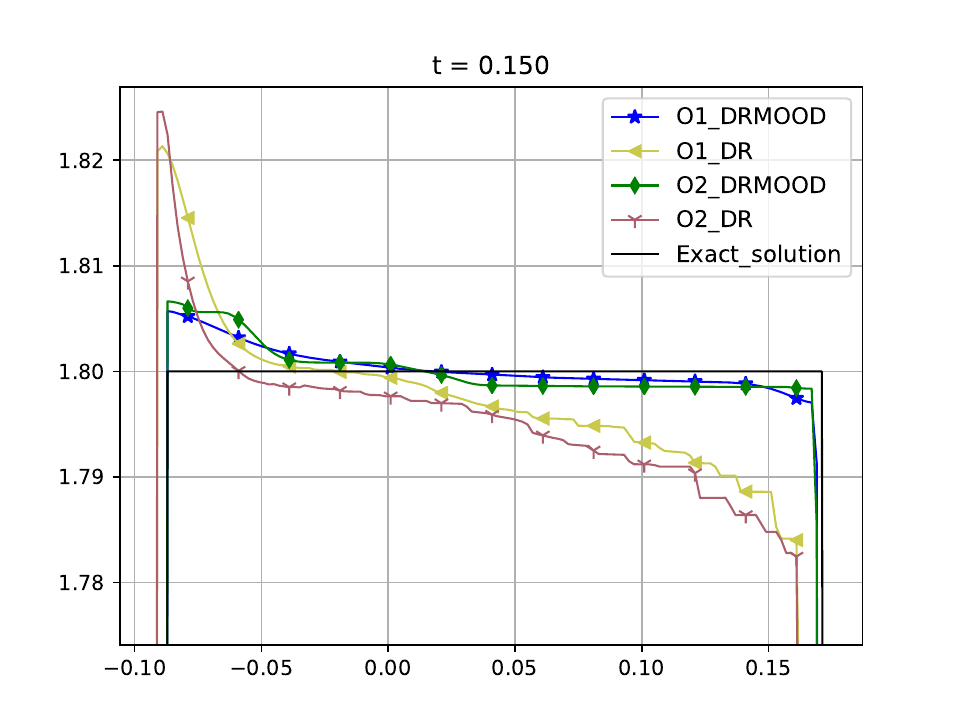}
				\caption*{Variable $h$: Zoom}
		\end{subfigure}
		\begin{subfigure}{0.5\textwidth}
			\includegraphics[width=1.1\linewidth]{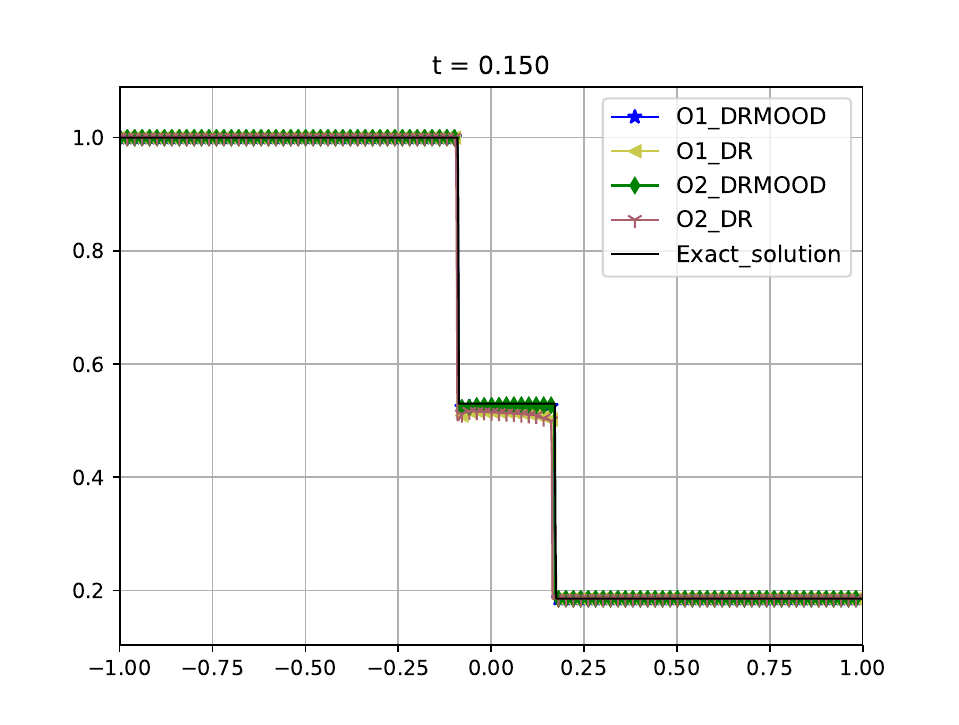}
			\caption{Variable $q$}
		\end{subfigure}
		\begin{subfigure}{0.5\textwidth}
				\includegraphics[width=1.1\linewidth]{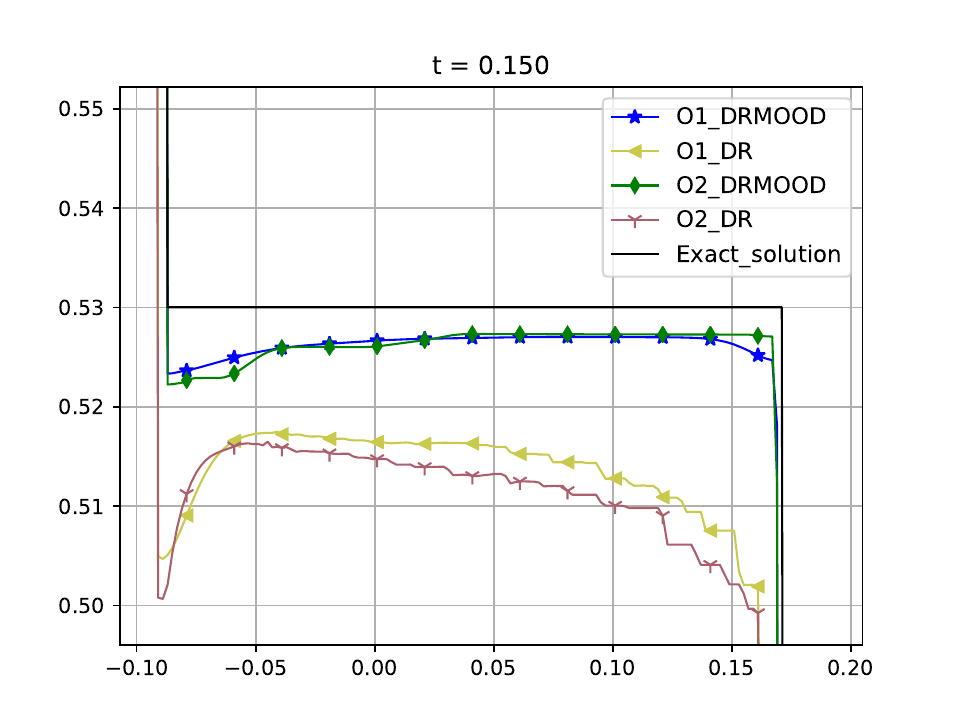}
				\caption*{Variable $q$: Zoom}
		\end{subfigure}
		\caption{Modified Shallow Water system. Test 2: Numerical solutions obtained using the first- and second-order in-cell discontinuous reconstruction methods with and without the MOOD approach at time $t=0.15$ with 1000 cells. Top: variable $h$ (left), zoom middle state (right). Down: variable $q$ (left), zoom middle state (right).}
		\label{fig:1DModifiedShallowWater_Test1_DR_vs_DRMOOD_1000_t015}
	\end{figure}

Let us consider the following initial condition taken from \cite{pimentel2021cell}
\begin{equation}\label{eq:1DS_Test2}
  \bu_0(x) = [h_0(x),q_{0}(x)]^T= \begin{cases}
     [1, 1]^T & \text{if $x<0$,}  \\
     [1.5,0.1855893974385]^T & \text{otherwise.}
\end{cases} 
\end{equation}
The solution of the Riemann problem consists of a 1-shock wave with negative speed and a 2-shock with positive speed  with intermediate state $\bu_* = [1.8, 0.530039370688997]^T$. Figure \ref{fig:1DModifiedShallowWater_Test1_1000_t015} compares the exact solution of the Riemann problem and the numerical results at time $t=0.15$ obtained with O$p$\_MOOD and O$p$\_DRMOOD, $p=1,2,3$, using a 1000-cell mesh. Notice that, in this case, the discontinuities are not isolated in the initial condition and thus the exact solution is not exactly captured by methods based on the in-cell discontinuous reconstruction, as it can be seen in the figure. Nevertheless, unlike the other methods, the numerical solutions converge to the right weak solution. 
In Figure \ref{fig:1DModifiedShallowWater_Test1_1000_vs_5000_t015} we compare the O$p$\_DRMOOD, $p=1,2,3$, using 1000- and 5000-cell meshes. We observe how all methods converge to the exact solution when the mesh is refined. In Figure \ref{fig:1DModifiedShallowWater_Test1_DR_vs_DRMOOD_1000_t015} we compare the exact solution and the numerical results at time $t=0.15$ obtained with O$p$\_DRMOOD and O$p$\_DR, $p=1,2,3$, using a 1000-cell mesh. We observe that all the methods seem to converge to the right solution but the ones based on the MOOD approach give better results. Moreover the CPU times are much better when using the DR.MOOD since the time steps are only restricted in the marked cells: see Table \ref{tab:tes2:runtimes}.

 \begin{table}[htb!]
    \center
	\begin{tabular}{ |c|c|c|c|c|c|c| }
		\hline
		Cells & O1\_DR & O1\_DRMOOD & Speedup & O2\_DR & O2\_DRMOOD & Speedup \\
		\hline
		500 & 4.468 & 0.46 & 9.73 & 37.08 & 0.45 & 81.43\\
		\hline
		1000 & 14.25 & 0.81 & 17.68 & 74.64 & 0.87 & 85.69\\
		\hline
		2000 & 44.79 & 2.10 & 21.32 & 173.71 & 2.27 & 76.45\\
		\hline
		4000 & 185.83 & 6.53 & 28.47 & 628.88 & 7.142 & 88.05\\
		\hline
	\end{tabular}
    \caption{Test 2: CPU times in (s) for different number of cells of the domain obtained for O$p$\_DRMOOD and O$p$\_DR, $p=1,2$ at time $t=0.5$ (average of 5 runs).}
    \label{tab:tes2:runtimes}
\end{table}

\subsubsection*{Test 3: right-moving 1-shock + right-moving 2-shock}\label{ssstest3}

\begin{figure}[htpb]
		\begin{subfigure}{0.5\textwidth}
			\includegraphics[width=1.1\linewidth]{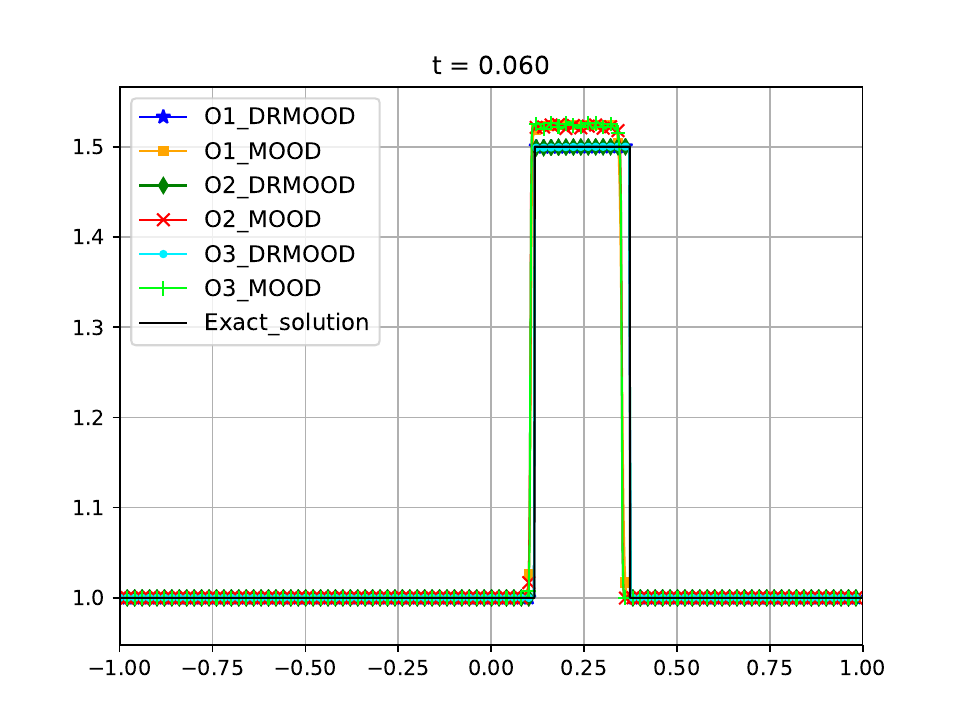}
			\caption{Variable $h$}
		\end{subfigure}
		\begin{subfigure}{0.5\textwidth}
				\includegraphics[width=1.1\linewidth]{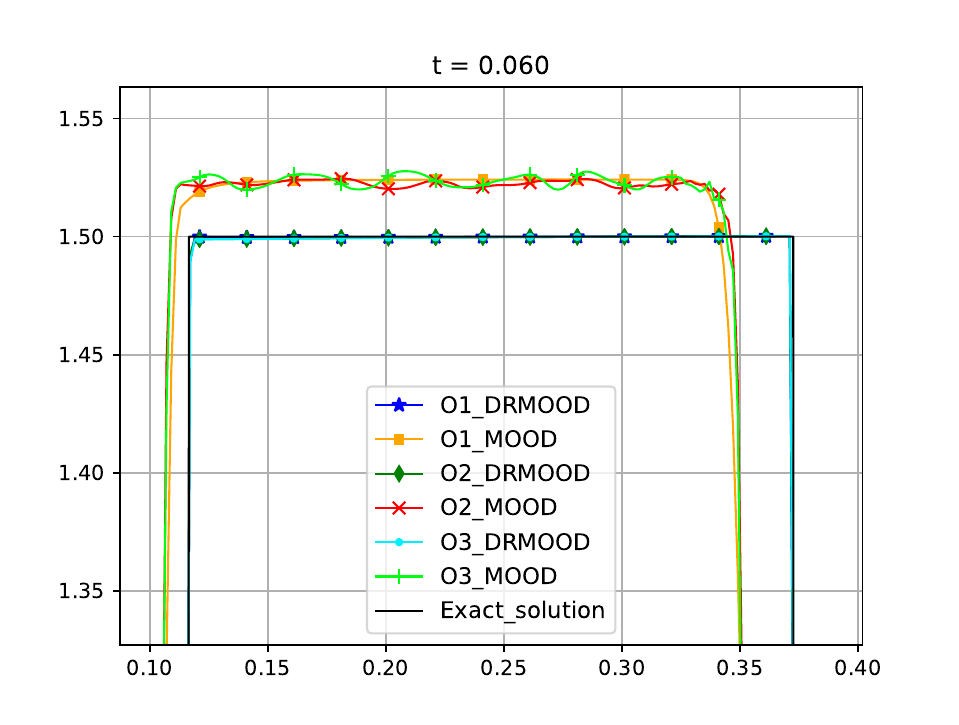}
				\caption*{Variable $h$: Zoom}
		\end{subfigure}
		\begin{subfigure}{0.5\textwidth}
			\includegraphics[width=1.1\linewidth]{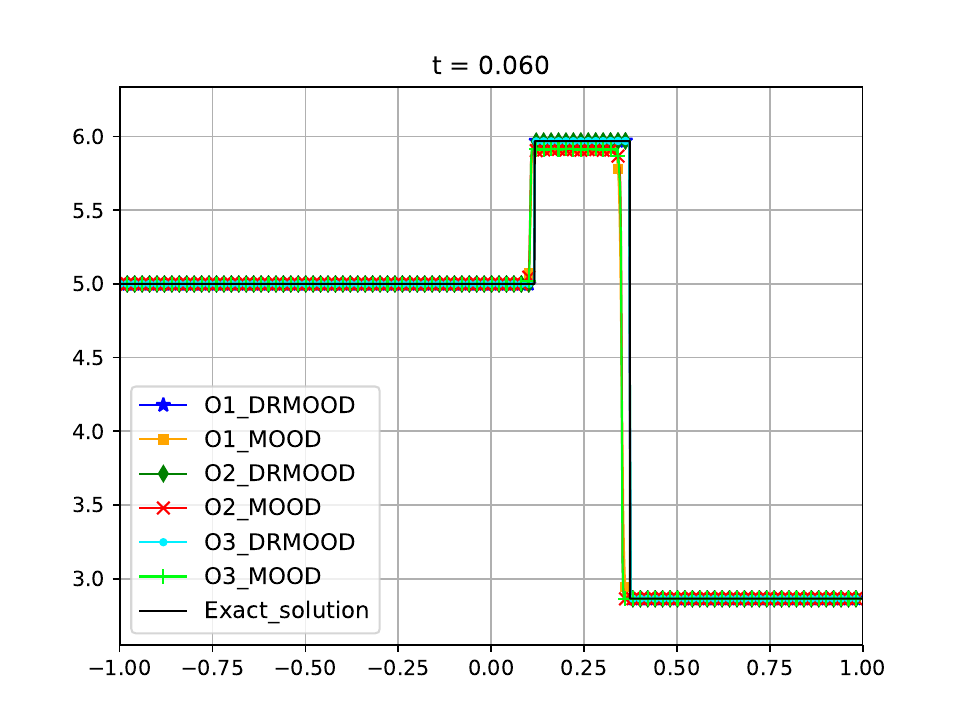}
			\caption{Variable $q$}
		\end{subfigure}
		\begin{subfigure}{0.5\textwidth}
				\includegraphics[width=1.1\linewidth]{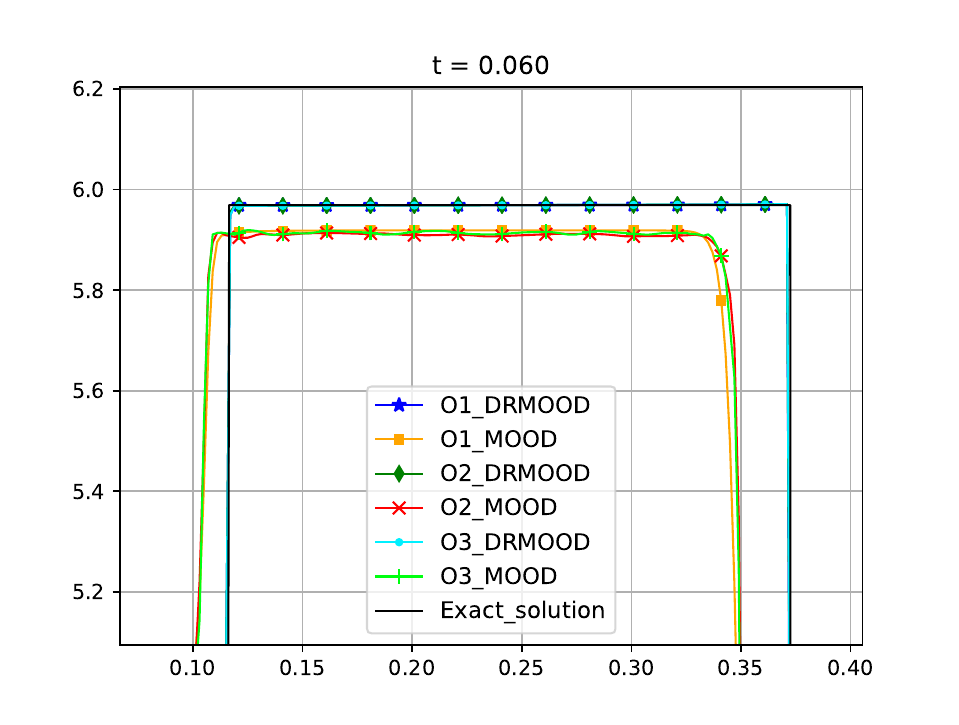}
				\caption*{Variable $q$: Zoom}
		\end{subfigure}
		\caption{Modified Shallow Water system. Test 3: Numerical solutions obtained using the MOOD approach with the first- , second- and third-order methods with and without in-cell discontinuous reconstruction based on the Roe matrix at time $t=0.06$ with 1000 cells. Top: variable $h$ (left), zoom middle state (right). Down: variable $q$ (left), zoom middle state (right).}
		\label{fig:1DModifiedShallowWater_Test2_1000_t006}
	\end{figure}

Let us consider the following initial condition taken from \cite{pimentel2021cell}

\begin{equation}\label{eq:1DS_Test3}
   \bu_0(x) = [h_0(x),q_{0}(x)]^T= \begin{cases}
     [1, 1]^T & \text{if $x<0$,}  \\
     [5,2.86423084288]^T & \text{otherwise.}
\end{cases}
\end{equation}
The solution of the Riemann problem consists of a 1-shock and a 2-shock waves with positive speed and intermediate state $\bu_* = [1.5, 5.96906891076]^T$. In Figure \ref{fig:1DModifiedShallowWater_Test2_1000_t006} we show the exact solution of the Riemann problem and the numerical results at time $t=0.06$ obtained with O$p$\_MOOD and O$p$\_DRMOOD, $p=1,2,3$, using a 1000-cell mesh. Observe that, again, the methods with the in-cell discontinuous reconstruction are able to capture with high accuracy the two shocks while the other ones are not converging to the right solution.

 \subsubsection*{Test 4: right-moving 1-rarefaction + right-moving 2-shock}\label{ssstest4}

 We consider the following initial condition

 \begin{equation}\label{eq:1DS_Test4}
  \bu_0(x) = [h_0(x),q_{0}(x)]^T= \begin{cases}
     [0.6, 0.3]^T & \text{if $x<0$,}  \\
     [1.5,0.7]^T & \text{otherwise.}
\end{cases} 
\end{equation}

\begin{figure}[htpb]
		\begin{subfigure}{0.5\textwidth}
			\includegraphics[width=1.1\linewidth]{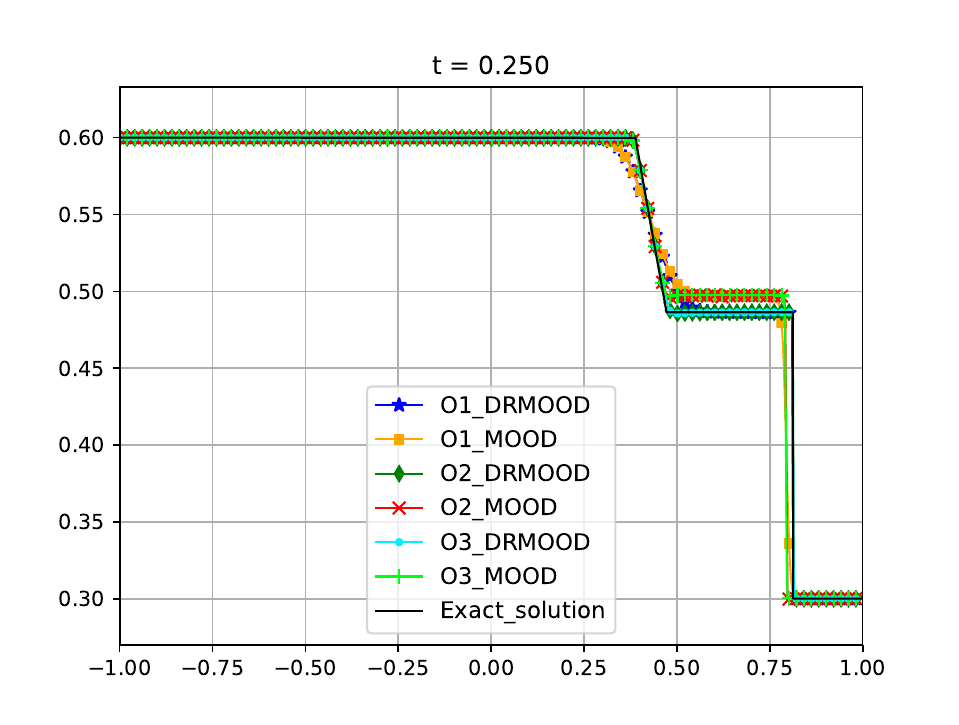}
			\caption{Variable $h$}
		\end{subfigure}
		\begin{subfigure}{0.5\textwidth}
				\includegraphics[width=1.1\linewidth]{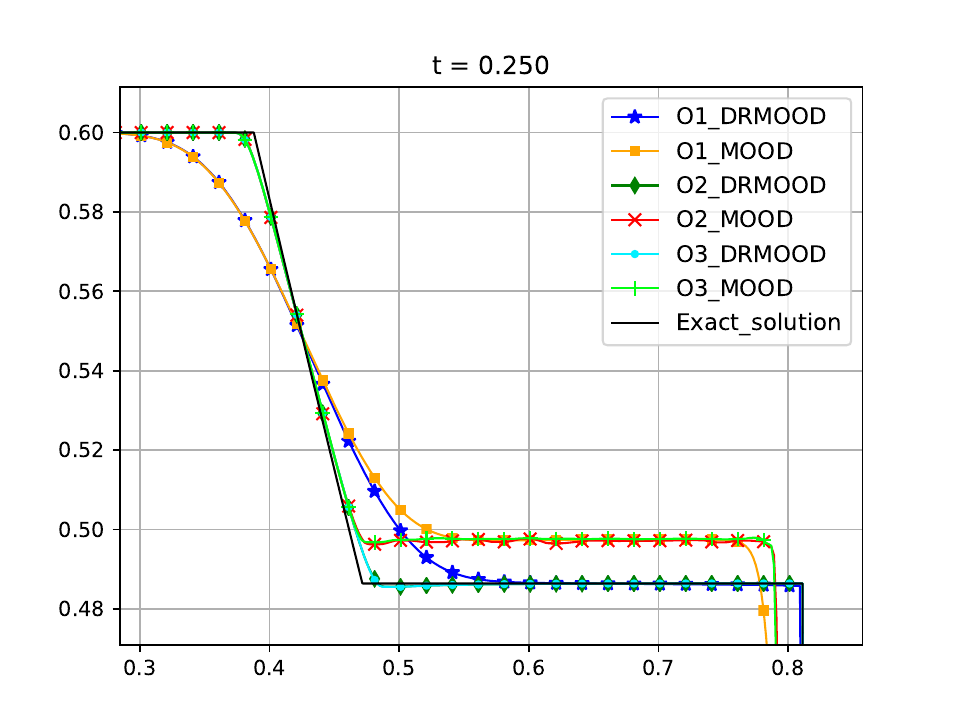}
				\caption*{Variable $h$: Zoom}
		\end{subfigure}
		\begin{subfigure}{0.5\textwidth}
			\includegraphics[width=1.1\linewidth]{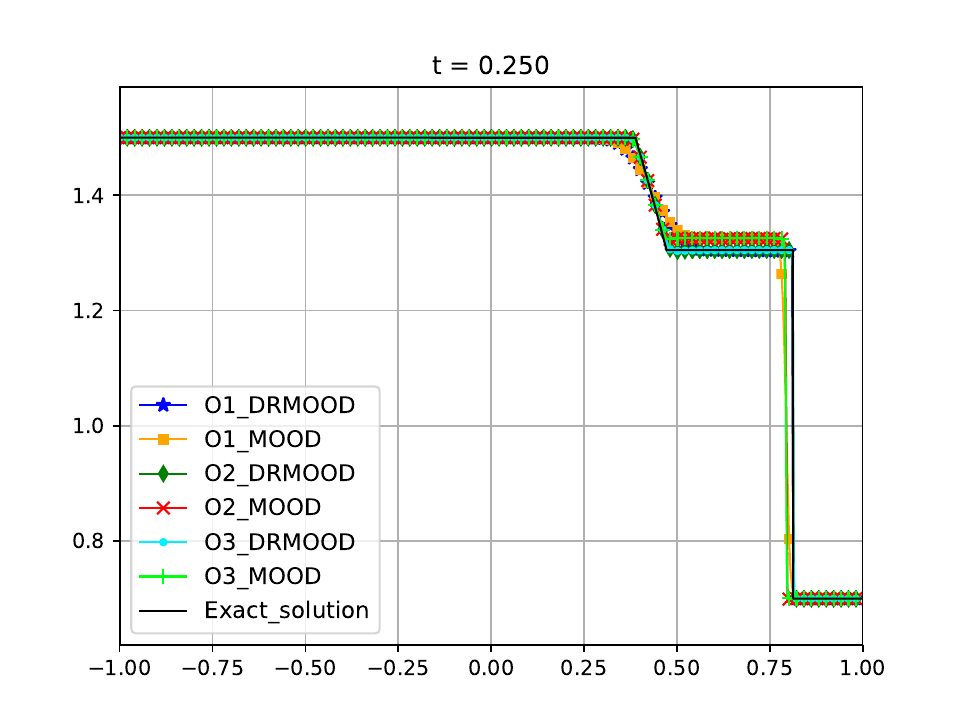}
			\caption{Variable $q$}
		\end{subfigure}
		\begin{subfigure}{0.5\textwidth}
				\includegraphics[width=1.1\linewidth]{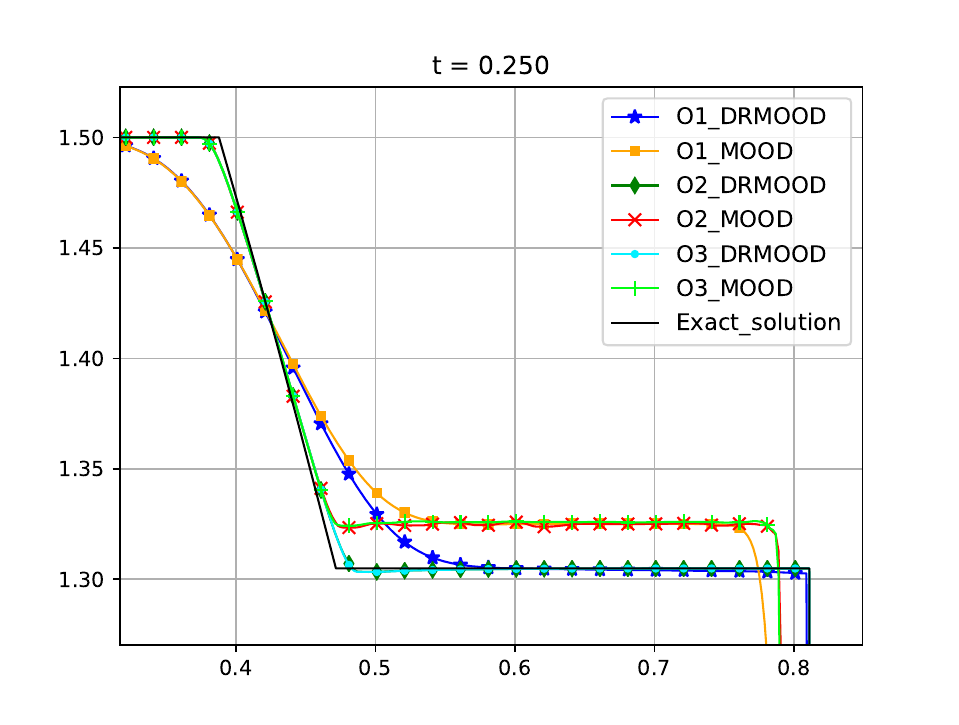}
				\caption*{Variable $q$: Zoom}
		\end{subfigure}
		\caption{Modified Shallow Water system. Test 4: Numerical solutions obtained using the MOOD approach with the first- , second- and third-order methods with and without in-cell discontinuous reconstruction based on the Roe matrix at time $t=0.25$ with 1000 cells. Top: variable $h$ (left), zoom middle state (right). Down: variable $q$ (left), zoom middle state (right).}
		\label{fig:1DModifiedShallowWater_Test3_1000_t025}
	\end{figure}

The solution of this Riemann problem consists of a 1-rarefaction whose head and tail speeds are positive and a 2-shock with positive speed too. The intermediate state is in this case $$\bu_* = [0.486393398885, 1.304922681678]^T.$$ 
Figure \ref{fig:1DModifiedShallowWater_Test3_1000_t025} compares the exact solution of the Riemann problem and the numerical results at time $t=0.25$ obtained with O$p$\_MOOD and O$p$\_DRMOOD, $p=1,2,3$, using a 1000-cell mesh. Again the methods with the in-cell discontinuous reconstructions are able to obtain the right convergence while the ones without it are not. We also observe a big improvement when capturing the rarefaction for the second- and third-order methods as expected. The difference between the second- and third-order schemes seems to be small due to the 2-shock in the solution. The following test is devoted to prove the second- and third-order accuracy of the methods in the smooth regions. 

\subsubsection*{Test 5: Order of accuracy}\label{ssstest5}
This test is devoted to show the order of accuracy of the methods. The initial condition consists on a constant state with a perturbation in the $h$ variable:
 \begin{equation}\label{eq:1DS_Test5}
  \bu_0(x) = [h_0(x) + \delta(x),q_{0}(x)]^T,
\end{equation}
where $h_0(x) = 1$, $q_{0}(x) = 0.5$ and
$$\delta(x) = 0.01\,e^{-100\,x^2}.$$

We use 100-, 200-, 400-, 800-, 1600- and 3200-cell uniform meshes in order to compute the errors in the $L^{1}$ norm and check the order of the first-, second- and third-order DR.MOOD schemes. The reference solution has been computed with the third-order DR.MOOD method using a 12800-point uniform grid. In Table \ref{tab:Error_TestOrder} we show the results at time $t=0.5$ from which we conclude that the expected order of accuracy is obtained in the three cases.

\begin{table}[!ht]
  	\centering
  	\begin{tabular}{|c|c|c|c|c|}
  	    \hline
         \multicolumn{5}{|c|}{First-order} \\
 \hline
  		Number of cells & $||\Delta h||_1$ & Order & $||\Delta q||_1$ & Order \\  
  		\hline 
   100 & 9.59e-04 & - & 3.00e-04 & - \\
   200 & 6.24e-04 & 0.62 & 1.92e-04 & 0.65 \\
   400 & 3.72e-04 & 0.75 & 1.12e-04 & 0.77 \\
   800 & 2.07e-04 & 0.84 & 6.18e-05 & 0.86 \\
  1600 & 1.10e-04 & 0.91 & 3.26e-05 & 0.92 \\
  3200 & 5.69e-05 & 0.95 & 1.68e-05 & 0.96 \\
  		\hline 
  		\hline
         \multicolumn{5}{|c|}{Second-order} \\
 \hline
  		Number of cells & $||\Delta h||_1$ & Order & $||\Delta q||_1$ & Order \\  
  		\hline 
   100 & 8.67e-05 & - & 2.44e-05 & - \\
   200 & 1.32e-05 & 2.72 & 3.62e-06 & 2.75 \\
   400 & 1.80e-06 & 2.87 & 4.87e-07 & 2.89 \\
   800 & 2.71e-07 & 2.73 & 7.08e-08 & 2.78 \\
  1600 & 5.09e-08 & 2.41 & 1.25e-08 & 2.50 \\
  3200 & 1.14e-08 & 2.16 & 2.63e-09 & 2.24 \\
  		\hline
  		\hline
         \multicolumn{5}{|c|}{Third-order} \\
 \hline
  		Number of cells & $||\Delta h||_1$ & Order & $||\Delta q||_1$ & Order \\  
  		\hline 
   100 & 6.24e-05 & - & 1.92e-05 & - \\
   200 & 8.94e-06 & 2.80 & 2.72e-06 & 2.82 \\
   400 & 1.15e-06 & 2.96 & 3.48e-07 & 2.97 \\
   800 & 1.44e-07 & 2.99 & 4.37e-08 & 2.99 \\
  1600 & 1.80e-08 & 3.00 & 5.46e-09 & 3.00 \\
  3200 & 2.22e-09 & 3.02 & 6.73e-10 & 3.02 \\
  		\hline
  	\end{tabular} 
	  	\caption{Test 5: order of accuracy for the first-, second- and third-order DR.MOOD scheme: $L^{1}$ errors $||\Delta \cdot||_1$ at time $t=0.5$.}

  	\label{tab:Error_TestOrder}
\end{table}
	
\subsection{Two-layer shallow water system}
We consider the homogeneous two-layer 1-D shallow water system (see \cite{castro2001q}):

\begin{equation}\label{eq:bicapa}
\left\{\begin{array}{l}
\displaystyle \frac{\partial h_{1}}{\partial t} + \frac{\partial q_{1}}{\partial x} = 0,\\
\displaystyle \frac{\partial q_{1}}{\partial t} + \frac{\partial}{\partial x}\left(\frac{q_{1}^{2}}{h_{1}} + \frac{1}{2}gh_{1}^{2}\right) = -gh_{1}\frac{\partial h_{2}}{\partial x},\\
\displaystyle \frac{\partial h_{2}}{\partial t} + \frac{\partial q_{2}}{\partial x} = 0,\\
\displaystyle \frac{\partial q_{2}}{\partial t} + \frac{\partial}{\partial x}\left(\frac{q_{2}^{2}}{h_{2}} + \frac{1}{2}gh_{2}^{2}\right) = -\frac{\rho_{1}}{\rho_{2}}gh_{2}\frac{\partial h_{1}}{\partial x}.
\end{array}\right.
\end{equation}

Index 1 refers to the upper layer while index 2 refers to the lower layer. This system uses the following notation:

\begin{itemize}
	\item $h_{i}=h_{i}(x,t)\geq 0$ is the thickness of the $i$-th layer at the section of coordinate $x$ at time t.
	\item $q_{i}=q_{i}(x,t)$ is the discharge of the $i$-th layer at the section of coordinate $x$ at time t.
	\item $g$ is the intensity of the gravitational field.
	\item $\rho_{i}$ refers to the constant density of the $i$-th layer.
\end{itemize}

The bottom is assumed to be flat. System (\ref{eq:bicapa}) can be rewritten in the form:
$$\partial_{t}\bu + \mathcal{A}(\bu)\partial_{x}\bu=0,$$
with
$$\bu= \left(\begin{array}{c}
h_{1}\\
h_{1}u_{1}\\
h_{2}\\
h_{2}u_{2}
\end{array}\right) \quad \mathcal{A}(\bu) = \left(\begin{array}{cccc}
0 & 1 & 0 & 0\\
gh_{1}-u_{1}^{2} & 2u_{1} & gh_{1} & 0 \\
0 & 0 & 0 & 1 \\
rgh_{2} & 0 & gh_{2}-u_{2}^{2} & 2u_{2}
\end{array}\right),
$$
where $r = \rho_{1}/\rho_{2}$ is the ratio of the layer densities and $u_{i}=q_{i}/h_{i}$, $i=1,2$ are the depth-averaged velocities.\\
We consider the Roe matrix $\mathcal{A}_{\Phi}$ of the system based on the family of straight segments
$$\Phi(s; \bu_{l}, \bu_{r}) = \bu_{l} + s(\bu_{r}-\bu_{l})$$
described in  \cite{pares2004well}. In this case we have 
$$\mathcal{A}_{\Phi}(\bu_{l}, \bu_{r}) = \mathcal{A}(\bar{\bu}),$$
where $\bar{\bu}=(\bar{h}_{1},\bar{h}_{1}\bar{u}_{1},\bar{h}_{2},\bar{h}_{2}\bar{u}_{2})^{T}$, being
$$\bar{h}_{p} = \frac{h_{p,l} + h_{p,r}}{2}, \quad \bar{u}_{p} = \frac{\sqrt{h_{p,l}}u_{p,l} + \sqrt{h_{p,r}}u_{p,r}}{\sqrt{h_{p,l}} + \sqrt{h_{p,r}}}, \quad p=1,2.$$

\subsubsection{Cell-marking criterion}

The variable $h_2$ is used to apply  the (DMP) \eqref{DMP} criterion for the second- and third-order schemes or the (LSJ) criterion \eqref{LSJ_criterion} for the first-order one. 

\subsubsection{Numerical tests}
Different numerical tests are presented to show the promising results given by the DR.MOOD methods. We consider a CFL number of 0.5, $g=9.81$, $r=\frac{\rho_1}{\rho_2} = 0.98$ and the domain $[0,3]$ in all tests.

\subsubsection*{Test 1: Isolated exterior shock}

Let us consider the following initial condition

\begin{eqnarray*}
& & \bu_0(x) = [h_{1,0}(x),q_{1,0}(x), h_{2,0}(x), q_{2,0}(x)]^T \\ & = & \begin{cases}
     [0.3701721263155, -0.186780084971819, 1.592710639376, 0.17351375487357]^T & \text{if $x<0.5$,}  \\
     [0.3612458594873,-0.22511224195307, 1.542922188541, -0.049717735973576]^T & \text{otherwise.}
\end{cases}
\end{eqnarray*}

\begin{figure}[htpb]

\begin{subfigure}{0.5\textwidth}
			\includegraphics[width=1.1\linewidth]{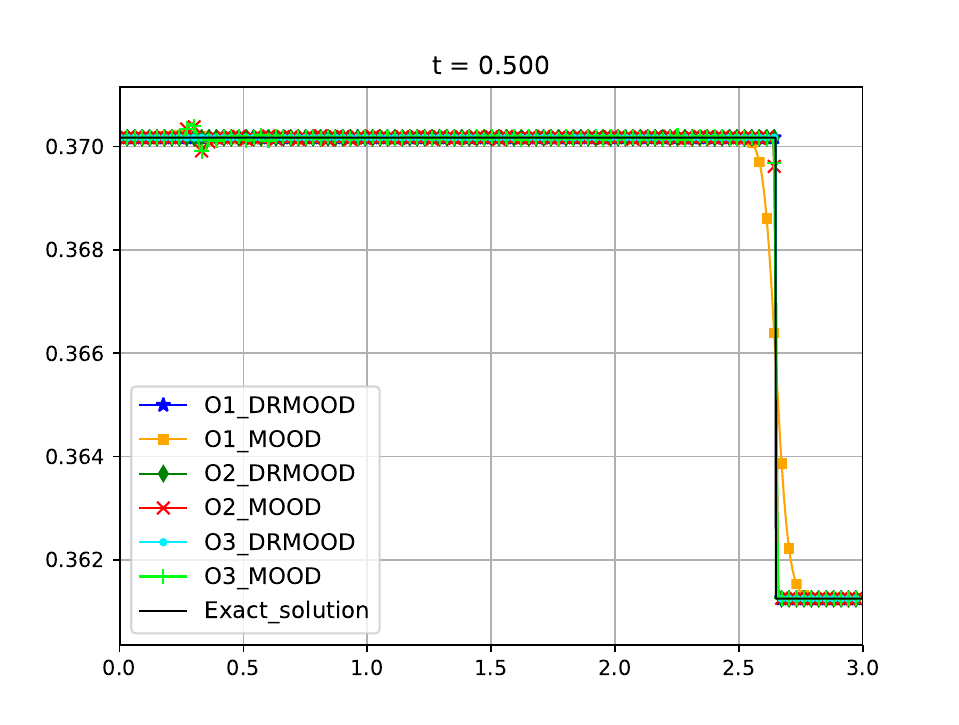}
			\caption{Variable $h_1$}
		\end{subfigure}
		\begin{subfigure}{0.5\textwidth}
				\includegraphics[width=1.1\linewidth]{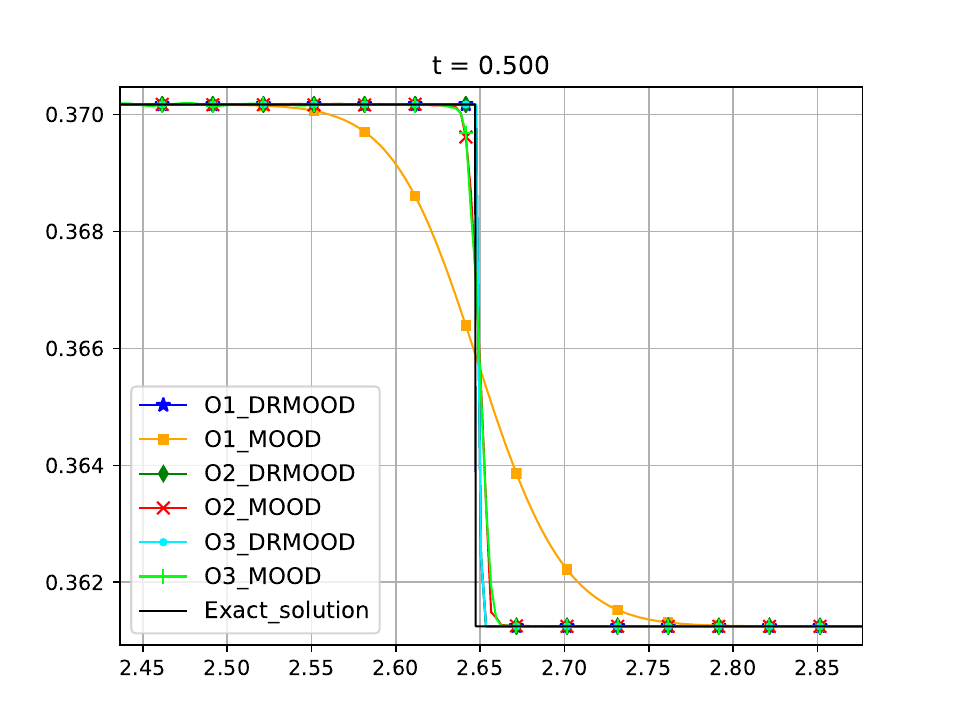}
				\caption*{Variable $h_1$: Zoom}
		\end{subfigure}
		\begin{subfigure}{0.5\textwidth}
			\includegraphics[width=1.1\linewidth]{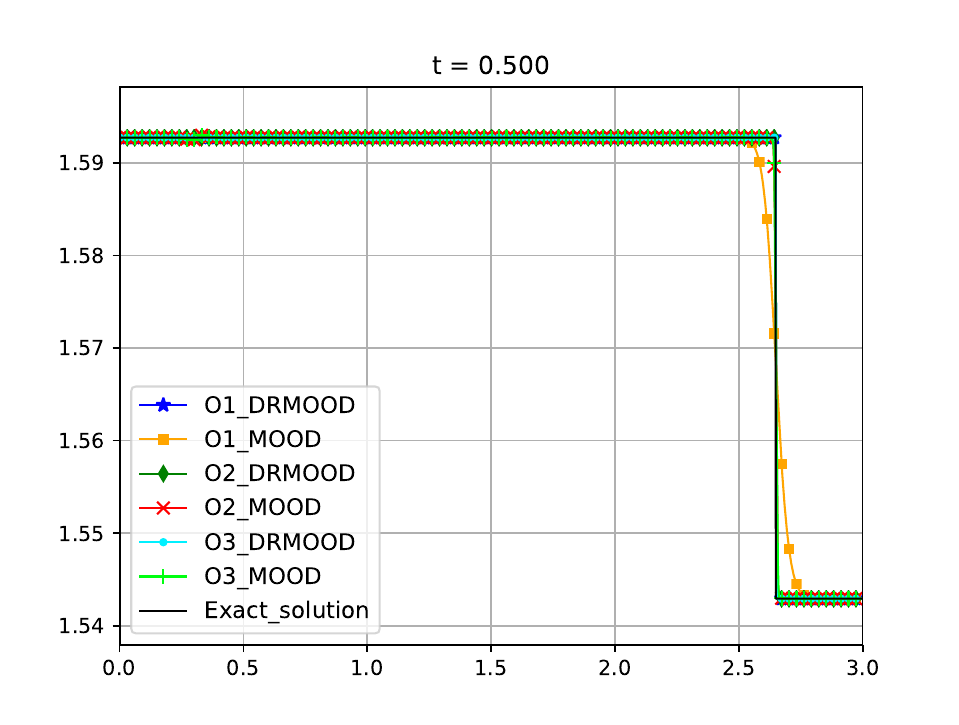}
			\caption{Variable $h_2$}
		\end{subfigure}
		\begin{subfigure}{0.5\textwidth}
				\includegraphics[width=1.1\linewidth]{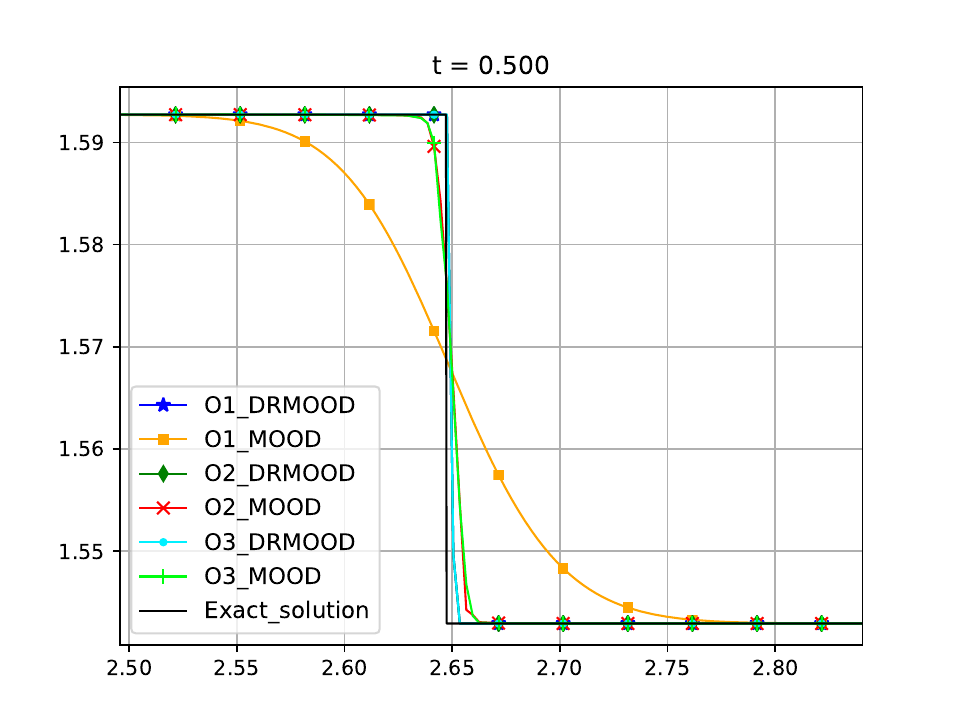}
				\caption*{Variable $h_2$: Zoom}
		\end{subfigure}
		\caption{Two-layer shallow water system. Test 1: Numerical solutions obtained using the MOOD approach with the first- , second- and third-order methods with and without in-cell discontinuous reconstruction based on the Roe matrix at time $t=0.5$ with 1000 cells. Top: variable $h_1$ (left), zoom (right). Down: variable $h_2$ (left), zoom (right).}
		\label{TwoLayerShallowWater_Test4Shock_1000_t05}
	\end{figure}

The solution of the Riemann problem consists of a 4-shock wave joining the left and right states. Figure \ref{TwoLayerShallowWater_Test4Shock_1000_t05} compares the exact solution of the Riemann problem and the numerical results at time $t=0.5$ obtained with O$p$\_MOOD and O$p$\_DRMOOD, $p=1,2,3$, using a 1000-cell mesh. In this case, all the numerical solutions are close to the exact one, although  the ones using the in-cell discontinuous reconstruction are more accurate. The fact that standard methods obtain the solution with small errors is known for  shocks related to the lowest and highest eigenvalues (that are called external shocks): see \cite{Castro2008}, where the numerical Hugoniot curves given by Roe and Rusanov schemes are compared to the exact Hugoniot curves.

\subsubsection*{Test 2: Isolated internal shock}

Let us consider the following initial condition
\begin{eqnarray*}
\bu_0(x) & = & [h_{1,0}(x),q_{1,0}(x), h_{2,0}(x), q_{2,0}(x)]^T \\ 
& = & \begin{cases}
     [1.164817, -0.0497756, 0.8134379, 0.0391596]^T & \text{if $x<0.5$,}  \\
     [0.370172126315573,-0.18678, 1.59271063937673, 0.173514]^T & \text{otherwise.}
\end{cases}
\end{eqnarray*}

\begin{figure}[htpb]

\begin{subfigure}{0.5\textwidth}
			\includegraphics[width=1.1\linewidth]{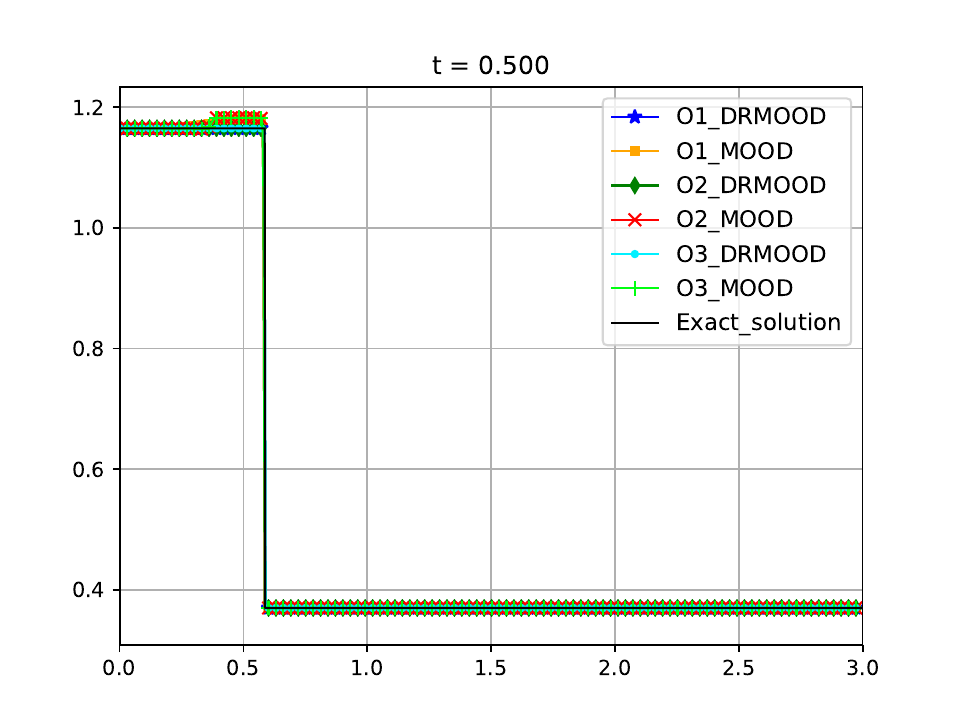}
			\caption{Variable $h_1$}
		\end{subfigure}
		\begin{subfigure}{0.5\textwidth}
				\includegraphics[width=1.1\linewidth]{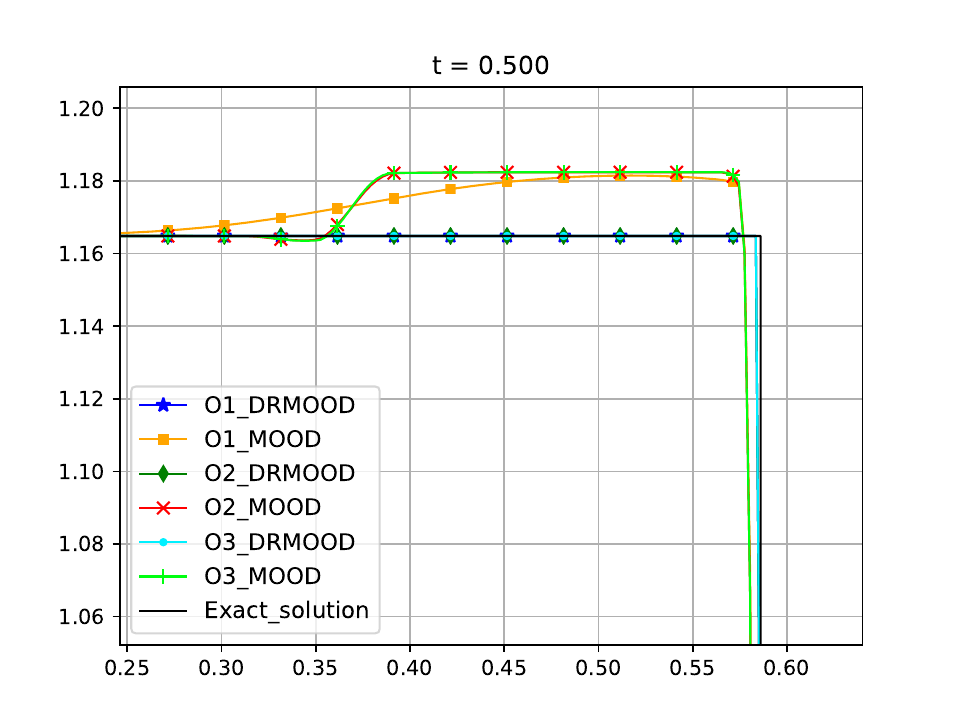}
				\caption*{Variable $h_1$: Zoom}
		\end{subfigure}
		\begin{subfigure}{0.5\textwidth}
			\includegraphics[width=1.1\linewidth]{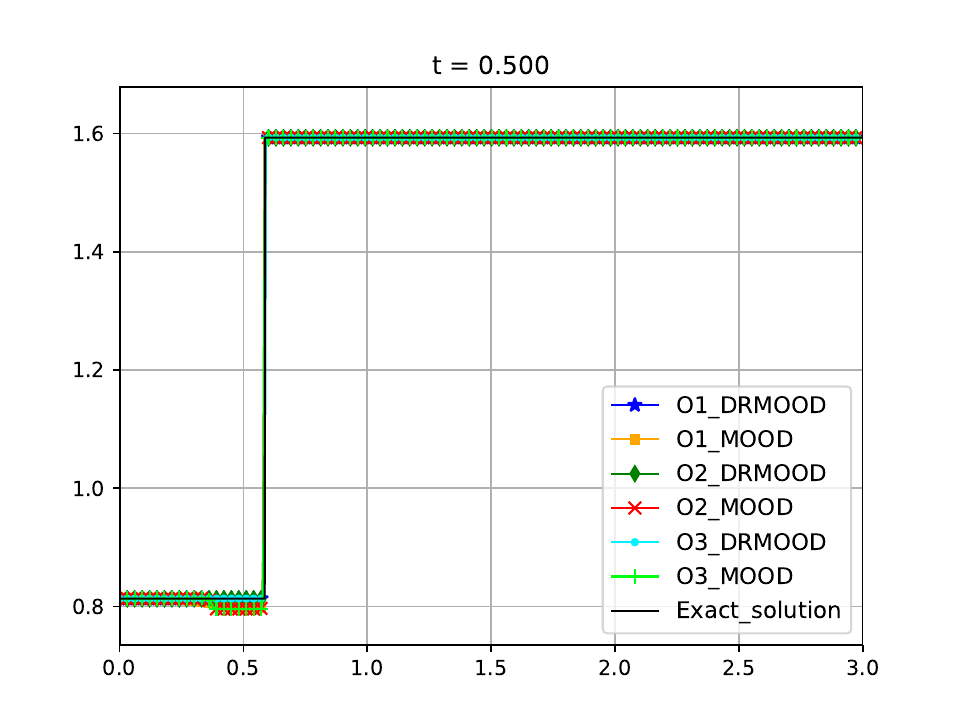}
			\caption{Variable $h_2$}
		\end{subfigure}
		\begin{subfigure}{0.5\textwidth}
				\includegraphics[width=1.1\linewidth]{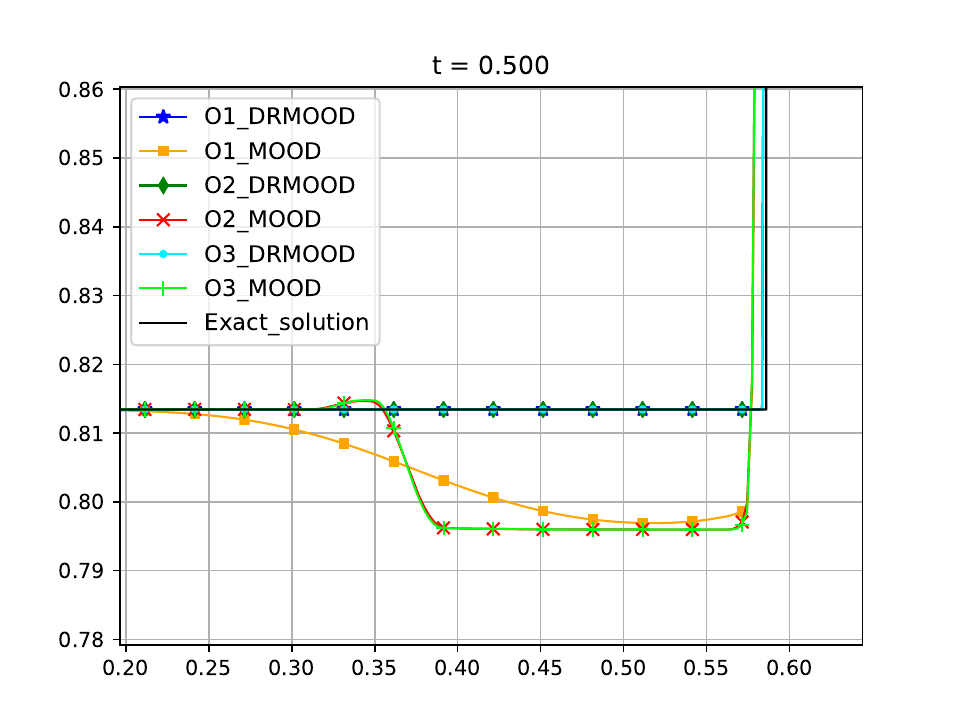}
				\caption*{Variable $h_2$: Zoom}
		\end{subfigure}
		\caption{Two-layer shallow water system. Test 2: Numerical solutions obtained using the MOOD approach with the first- , second- and third-order methods with and without in-cell discontinuous reconstruction based on the Roe matrix at time $t=0.5$ with 1000 cells. Top: variable $h_1$ (left), zoom (right). Down: variable $h_2$ (left), zoom (right).}
		\label{TwoLayerShallowWater_Test3Shock_1000_t05}
	\end{figure}

The solution of the Riemann problem consists of a 3-shock wave joining the left and right states.
Figure \ref{TwoLayerShallowWater_Test3Shock_1000_t05} shows the exact solution of the Riemann problem and the numerical results at time $t=0.5$ obtained with O$p$\_MOOD and O$p$\_DRMOOD, $p=1,2,3$, using a 1000-cell mesh. We observe that the methods without the in-cell discontinuous reconstructions are not able to capture the 3-shock wave correctly while the other ones are able to capture it exactly. The fact that standard methods give bigger errors for shock waves related to the intermediate eigenvalues (that are called internal shocks) is also known: see the comparison of the numerical and exact Hugoniot curves in \cite{Castro2008}.

 \subsubsection*{Test 3: Two shocks}

Let us consider the following initial condition

\begin{eqnarray*}
& & \bu_0(x) = [h_{1,0}(x),q_{1,0}(x), h_{2,0}(x), q_{2,0}(x)]^T \\
& = & \begin{cases}
     [1.164817, -0.0497756, 0.8134379, 0.0391596]^T & \text{if $x<0.5$,}  \\
     [0.3612458594874,-0.22511224195308, 1.542922188541, -0.040293287871649]^T & \text{otherwise.}
\end{cases}
\end{eqnarray*}

\begin{figure}[htpb]

\begin{subfigure}{0.33\textwidth}
			\includegraphics[width=1.1\linewidth]{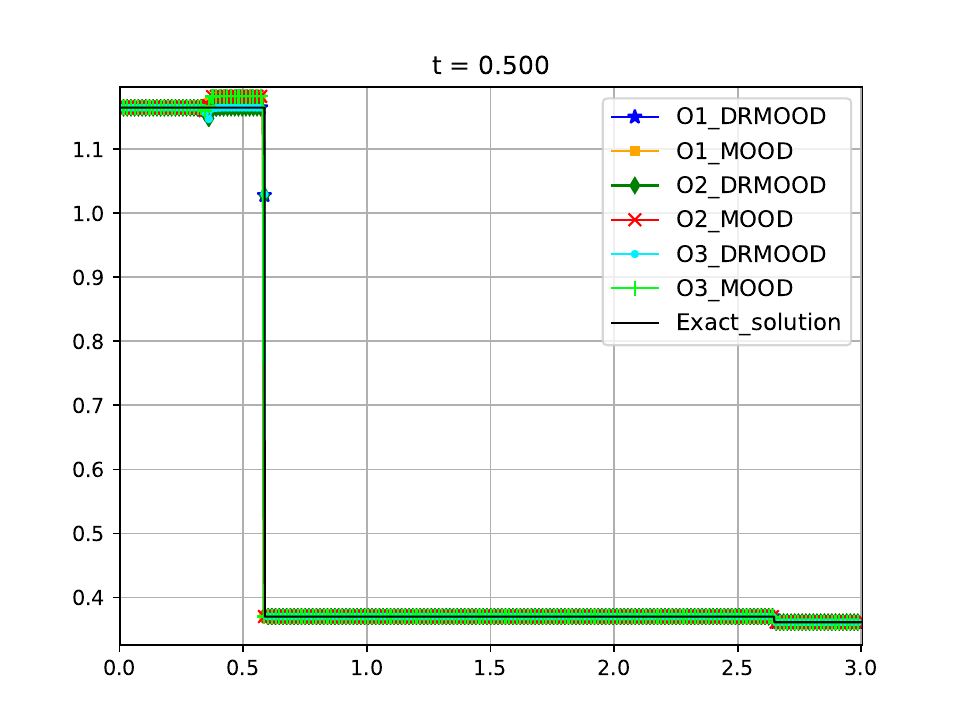}
			\caption{Variable $h_1$}
		\end{subfigure}
		\begin{subfigure}{0.33\textwidth}
				\includegraphics[width=1.1\linewidth]{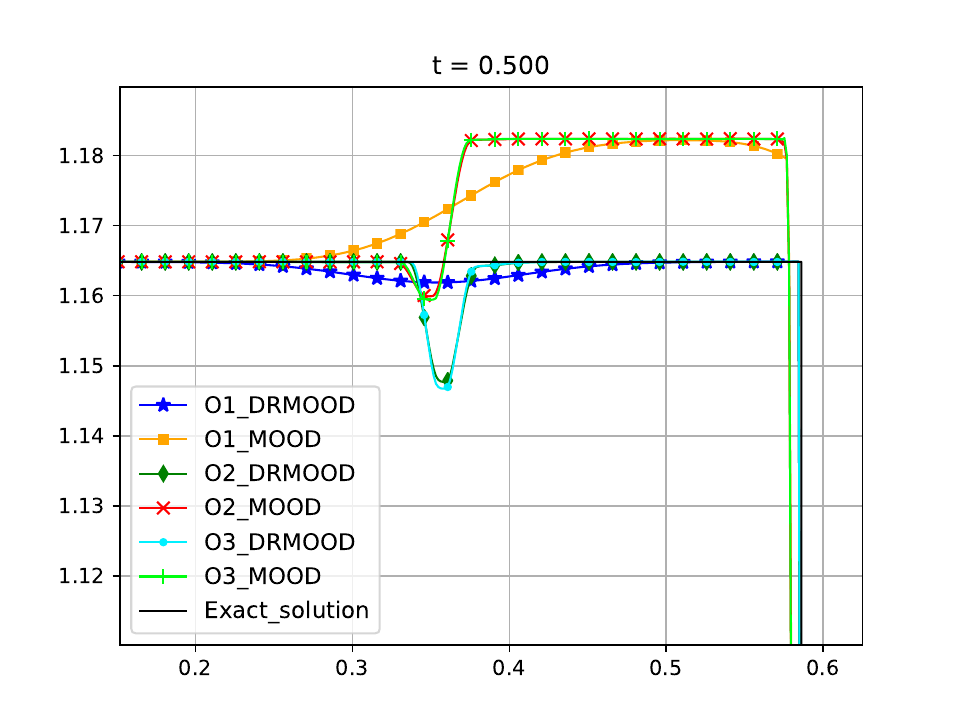}
				\caption*{Variable $h_1$: 3-shock zoom}
		\end{subfigure}
  \begin{subfigure}{0.33\textwidth}
				\includegraphics[width=1.1\linewidth]{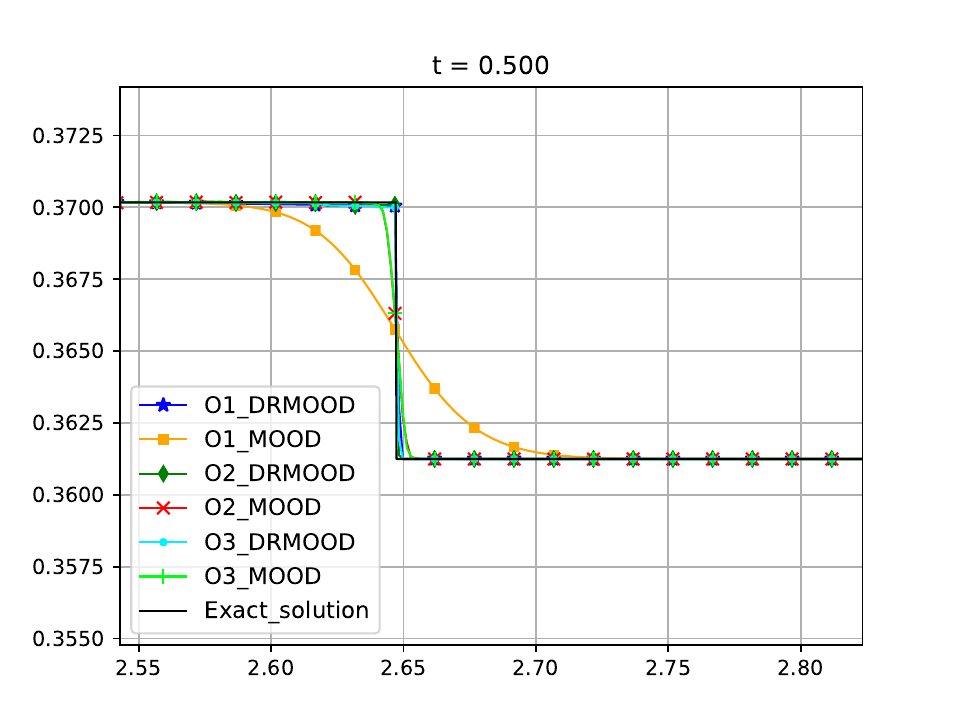}
				\caption*{Variable $h_1$: 4-shock zoom}
		\end{subfigure}
		\begin{subfigure}{0.33\textwidth}
			\includegraphics[width=1.1\linewidth]{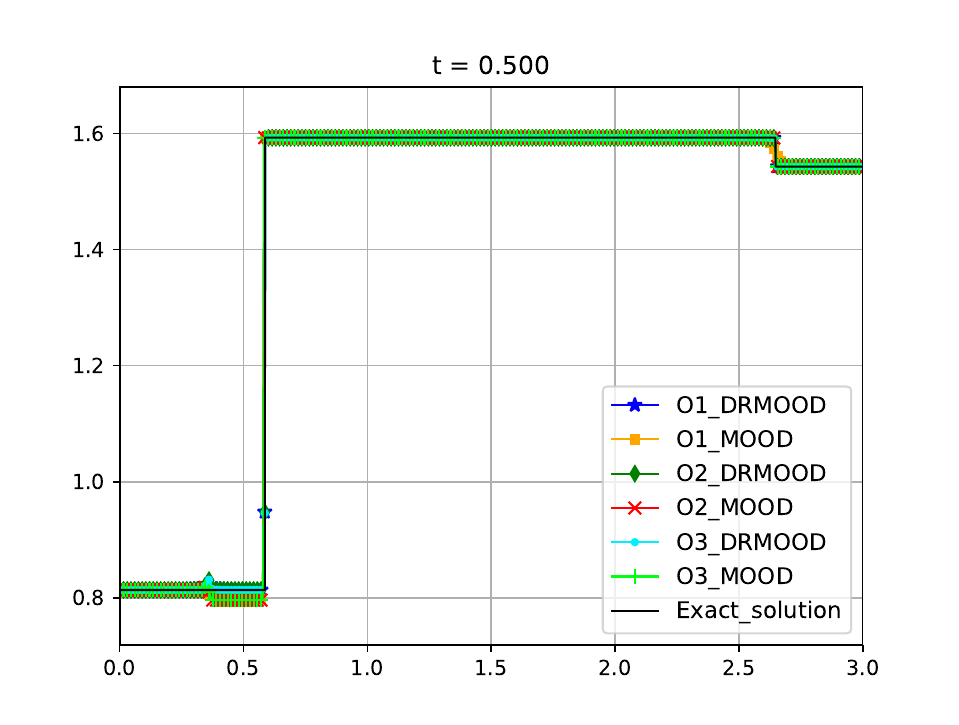}
			\caption{Variable $h_2$}
		\end{subfigure}
		\begin{subfigure}{0.33\textwidth}
				\includegraphics[width=1.1\linewidth]{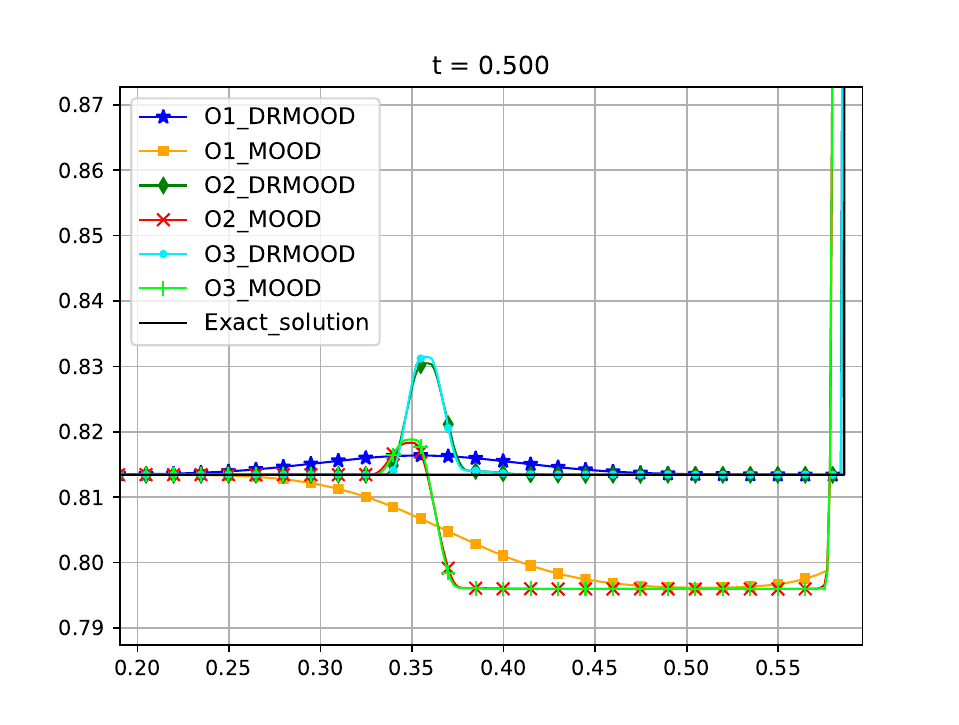}
				\caption*{Variable $h_2$: 3-shock zoom}
		\end{subfigure}
  \begin{subfigure}{0.33\textwidth}
				\includegraphics[width=1.1\linewidth]{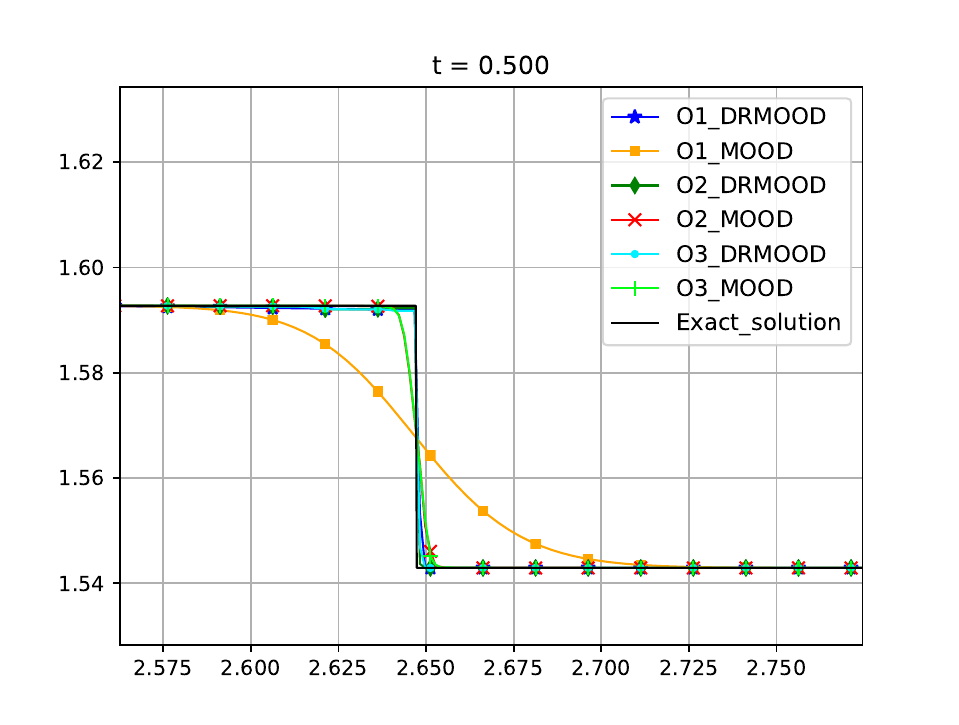}
				\caption*{Variable $h_2$: 4-shock zoom}
		\end{subfigure}
		\caption{Two-layer shallow water system. Test 3: Numerical solutions obtained using the MOOD approach with the first- , second- and third-order methods with and without in-cell discontinuous reconstruction based on the Roe matrix at time $t=0.5$ with 2000 cells. Top: variable $h_1$ (left), 3-shock zoom (middle), 4-shock zoom (right). Down: variable $h_2$ (left), 3-shock zoom (middle), 4-shock zoom (right).}
		\label{TwoLayerShallowWater_Test34Shock_2000_t05}
	\end{figure}

The solution of the Riemann problem consists of a 3-shock with negative speed and a 4-shock with positive speed. The intermediate state is given by 
$$\bu_* = [0.370172126315573, -0.18678008497181986 , 1.59271063937673, 0.1735137548735771]^T.$$
Figure \ref{TwoLayerShallowWater_Test34Shock_2000_t05} shows the exact solution of the Riemann problem and the numerical results at time $t=0.5$ obtained with O$p$\_MOOD and O$p$\_DRMOOD, $p=1,2,3$, using a 2000-cell mesh. As we can see on the right zooms, all methods are able to capture well the 4-shock although the ones with the in-cell discontinuous reconstructions are more accurate. In the case of the 3-shock wave the only methods that converge to the correct solution are those with the in-cell discontinuous reconstructions. A small oscillation appears in DR.MOOD methods close to the 3-shock that decresases as $\Delta x \rightarrow 0$: since the two shocks are not initially isolated, the solution is not exactly captured.

 \subsection{Conclusions}

 In this paper, a new family of path-conservative methods to deal numerically with discontinuous solutions of nonconservative hyperbolic systems is introduced. First the MOOD \cite{clain2011high} methodology is described for nonconservative systems, giving a detailed description of the correction terms added to be path-conservative. We combine this extension with the path-conservative in-cell discontinuous reconstruction methods introduced in \cite{chalons2019path}. The resulting scheme, DR.MOOD, is able to maintain the shock-capturing property for isolated shocks and it improves the results in the smooth regions of the solution. Several numerical tests are shown for the modified shallow water equations and the two-layer shallow water system. From these tests we observe that DR.MOOD methods converge to the right weak solutions. This new framework confirms what it was pointed out in \cite{Castro2008}: path-conservative methods converge to the correct solutions when the numerical viscosity close to discontinuities is correctly controlled. The main advantages of DR.MOOD strategy are:

 \begin{itemize}
     \item It is a general strategy: the only requirement to apply it that a  Roe matrix consistent with the selected family of paths has to be available.
     \item It captures exactly isolated shocks, improves the results in the smooth parts and it seems to converge to the correct solutions when non-isolated discontinuity appears.
     \item It leads to numerical methods that are more efficient than the high-order in-cell discontinuous reconstruction schemes introduced in \cite{pimentel2021cell}, since the time step is only reduced in the cells close to a shock.
 \end{itemize}
 
Some future works can be considered: application to other 1d nonconservative systems, analysis of the well-balanced property, analysis of the convergence, and extension to multidimensional problem.

\appendix
\section{DR.MOOD: correction terms in a boundary cell $j \in \mathcal{B}_{n}$}\label{app}
Let us consider  a cell $j\in \mathcal{B}_{n}$ such that $j+1\in \mathcal{O}_{n}$. 
The computation of $\bu_{j+1}^{n + 1}$ is done through $k$ intermediate steps of the in-cell reconstruction method:
\begin{eqnarray*}
    \bu_{j+1}^{n,i+1} & = & \bu_{j+1}^{n,i} - \frac{\Delta t_{n,i}}{\Delta x}\Bigl(\mathcal{D}_{j+\frac{3}{2}}^{-,i,DR} + \mathcal{D}_{j+\frac{1}{2}}^{+,i,DR}  \\
     & & \qquad \qquad  +  \dashint_{x_{j+ \frac{1}{2}}}^{x_{j+\frac{3}{2}}}\mathcal{A}(P_{j+1}^{DR,n,i}(x,t^{n}_i))
    \partial_x P_{j+1}^{DR,n,i}(x,t^{n}_i)\, dx\Bigr),\quad i = 0, \dots, k-1,
\end{eqnarray*}
where $P_{j+1}^{DR,n,i}$ represents the in-cell discontinuous reconstruction computed at the $(j+1)$th cell at time $t^{n}_i$ and
$$
\mathcal{D}_{j+\frac{1}{2}}^{+,i,DR}  = \mathcal{D}^{+,DR} \left( \bu_{j+\frac{1}{2}}^{-,HO}(t^{n}_i),
\bu_{j+\frac{1}{2}}^{+,LO}(t^{n}_i)\right),$$
where the following notation is used:
$$
\bu_{j+\frac{1}{2}}^{-,HO}(t^{n}_i) = P^n_j(x_{i+1/2}, t^{n}_i), \quad 
\bu_{j+\frac{1}{2}}^{+,LO}(t^{n}_i)  = P^{DR,n,i}_{j+1}(x_{i+1/2}, t^{n}_i),
$$
where $P^n_j(x,t)$ is the polynomial reconstruction corresponding to the high-order method computed in the $j$th cell at time $t^n$.
In the case of first-order DR.MOOD methods, in which such a polynomial is not needed, the first-degree interpolation polynomial
$$
P^n_j (x,t) = \bu_j^{n} + \frac{t- t^n}{\Delta t^{CFL}}\left( \bar\bu_j^{n +1} - \bu_j^n \right)$$
is considered.

Once these intermediate steps have been computed, the numerical solution is updated by taking $\bu_{j+1}^{n + 1} = \bu_{j+1}^{n,k}$. Summing up, the total contribution of the nonconservative products at the right of intercell $x_{j+1/2}$ is then:
$$
\sum_{i} \Delta t_{n,i}\mathcal{D}_{j+\frac{1}{2}}^{+,i,DR}.
$$
As for MOOD methods,  $\bu_j^{n +1}$ is corrected so that the left contribution of the nonconservative products is equal to
$$
\sum_{i} \Delta t_{n,i}\mathcal{D}_{j+\frac{1}{2}}^{-,i,DR}
$$
and so that the contributions related to all the jumps between $\bu_{j+\frac{1}{2}}^{-,HO}(t^{n,l}) $ and
$\bu_{j+\frac{1}{2}}^{+,LO}(t^{n}_i) $ are taken into account what, following the corrections in the MOOD method \eqref{mood_correction_LO}, gives
\begin{eqnarray}\label{dr_mood_correction_LO}
    \bu_{j}^{n+1}  =  \bu_{j}^{n} & - &  \frac{\Delta t_{n}^{CFL}}{\Delta x} \sum_{l=1}^{M_{t}} \beta_{l}^{t}\mathcal{D}_{j-\frac{1}{2}}^{+,HO}(t^{n,l})  - \sum_{i}\frac{\Delta t_{n,i}}{\Delta x}\mathcal{D}_{j+\frac{1}{2}}^{-,i,DR}  \\ \nonumber & - & \Delta  t_{n}^{CFL} \sum_{l=1}^{M_{t}} \sum_{m=1}^{M_{x}} \beta_{l}^{t} \beta_{m}^{x}\mathcal{A}(P_{j}^{n}(x_{j}^{m},t^{n,l}))
    \partial_x P_{j}^{n}(x_{j}^{m},t^{n,l}) \\   \nonumber & - &  \sum_{i}\frac{\Delta t_{n,i}}{\Delta x}\sum_{l=1}^{M_{t}} \beta_{l}^{t}\int_0^1\mathcal{A}(\Phi(s;\bu_{j+\frac{1}{2}}^{-,HO}(t^{n,l}), \bu_{j+\frac{1}{2}}^{+,DR}(t^{n}_i))\frac{\partial\Phi}{\partial
s}(s;\bu_{j+\frac{1}{2}}^{-,HO}(t^{n,l}), \bu_{j+\frac{1}{2}}^{+,DR}(t^{n}_i))\,ds
\end{eqnarray}
which can also be rewritten as follows: 
\begin{eqnarray}\label{dr_mood_correction_LO_2}
    \bu_{j}^{n+1}  =  \bu_{j}^{n} & - &  \frac{\Delta t_{n}^{CFL}}{\Delta x} \sum_{l=1}^{M_{t}} \beta_{l}^{t}\mathcal{D}_{j-\frac{1}{2}}^{+,HO}(t^{n,l})  - \sum_{i}\frac{\Delta t_{n,i}}{\Delta x} \mathcal{D}_{j+\frac{1}{2}}^{-,i,DR} \\ \nonumber & - &  \Delta  t_{n}^{CFL} \sum_{l=1}^{M_{t}} \sum_{m=1}^{M_{x}} \beta_{l}^{t} \beta_{m}^{x}\mathcal{A}(P_{j}^{n}(x_{j}^{m},t^{n,l}))
    \partial_x P_{j}^{n}(x_{j}^{m},t^{n,l}) \\ \nonumber & - &  \sum_{i}\frac{\Delta t_{n,i}}{\Delta x}\sum_{l=1}^{M_{t}} \beta_{l}^{t}\left(\mathcal{D}^{-}(\bu_{j+\frac{1}{2}}^{-,HO}(t^{n,l}), \bu_{j+\frac{1}{2}}^{+,DR}(t^{n}_i))) + \mathcal{D}^{+}(\bu_{j+\frac{1}{2}}^{-,HO}(t^{n,l}), \bu_{j+\frac{1}{2}}^{+,DR}(t^{n}_i))\right),
\end{eqnarray}
where the fluctuations $\mathcal{D}^{\pm}$ are taken as the Roe method \eqref{NCFR} since we selected the second strategy in Subsection \ref{ss:m_n}.

A similar expression is found for the correction in a cell $j\in \mathcal{B}_n$ such that $j-1\in \mathcal{O}_{n}$. An important remark is that again these corrections lead to a numerical method that is conservative if the system is conservative, i.e. if 
$\mathcal{A}(\bu)$ is the Jacobian of a flux function $F(\bu)$. The proof is similar to the one given in the MOOD procedure \eqref{proof_conservation}.


\clearpage

\bibliographystyle{abbrv}

\bibliography{references}

\end{document}